\documentclass[10pt]{amsart}
\usepackage[T1]{fontenc}
\usepackage{geometry}
\usepackage[latin1] {inputenc}
\usepackage{amsmath}
\usepackage{amsfonts, amssymb, textcomp}
\usepackage[colorlinks=flase, linkcolor=red,urlcolor=green, citecolor=blue]{hyperref}
\usepackage{subeqnarray}

\usepackage{latexsym}
\usepackage{fancyhdr}
\usepackage{longtable}
\usepackage{amsmath, amssymb}
\usepackage{graphicx}
\usepackage{enumerate}

\numberwithin{equation}{section}

\theoremstyle{plain}
\newtheorem{proposition}{Proposition}[section]
\newtheorem{theorem}[proposition]{Theorem}
\newtheorem{lemma}[proposition]{Lemma}
\newtheorem{corollary}[proposition]{Corollary}

\newtheorem{remark}[proposition]{Remark}

\newcommand{\RR}{\mathbb{R}}
\newcommand{\CC}{\mathbb{C}}
\newcommand{\NN}{\mathbb{N}}

\let\on=\operatorname

\newenvironment{demo}[1]{\par\smallskip\noindent{\bf #1.}}{\par\smallskip}

\makeatletter
\@namedef{subjclassname@2020}{%
  \textup{2020} Mathematics Subject Classification}
\makeatother

\newsavebox{\fmbox}
\newenvironment{fmpage}[1]
 {\begin{lrbox}{\fmbox}\begin{minipage}{#1}}
 {\end{minipage}\end{lrbox}\fbox{\usebox{\fmbox}}}

\title[Equality of ultradifferentiable classes by means of indices]
{Equality of ultradifferentiable classes by means of indices of mixed O-regular variation}

\author[J. Jim\'{e}nez-Garrido, J. Sanz, and G.~Schindl]{Javier Jim\'{e}nez-Garrido, Javier Sanz, and Gerhard Schindl}

\begin{document}

\begin{abstract}
We characterize the equality between ultradifferentiable function classes defined in terms of abstractly given weight matrices and in terms of the corresponding matrix of associated weight functions by using new growth indices. These indices, defined by means of weight sequences and (associated) weight functions, are extending the notion of O-regular variation to a mixed setting. Hence we are extending the known comparison results concerning classes defined in terms of a single weight sequence and of a single weight function and give also these statements an interpretation expressed in O-regular variation.
\end{abstract}

\keywords{Classes of ultradifferentiable functions; weight sequences, functions and matrices; growth indices; O-regular variation; mixed setting}
\subjclass[2020]{primary 46E10; secondary 26A12, 26A48, 46A13}
\date{\today}

\maketitle

\section{Introduction}\label{Introduction}
In the theory of ultradifferentiable function spaces there exist two classical approaches in order to control the growth of the derivatives of the functions belonging to such classes: Either one uses a weight sequence $M=(M_j)_j$ or a weight function $\omega:[0,+\infty)\rightarrow[0,+\infty)$. In both settings one requires several basic growth and regularity assumptions on $M$ and $\omega$ and one distinguishes between two types, the {\itshape Roumieu type spaces} $\mathcal{E}_{\{M\}}$ and $\mathcal{E}_{\{\omega\}}$, and the {\itshape Beurling type spaces} $\mathcal{E}_{(M)}$ and $\mathcal{E}_{(\omega)}$. In the following we write $\mathcal{E}_{[\star]}$ for all arising classes of ultradifferentiable functions if we mean either $\mathcal{E}_{\{\star\}}$ or $\mathcal{E}_{(\star)}$, but not mixing the cases. Similarly this is done for all arising conditions.\vspace{6pt}

In \cite{BonetMeiseMelikhov07} both methods have been compared: A characterization has been given when $\mathcal{E}_{[M]}=\mathcal{E}_{[\omega]}$ is valid (as locally convex vector spaces) and it has been shown that in general both approaches are mutually distinct.\vspace{6pt}

Motivated by the results from \cite{BonetMeiseMelikhov07}, in \cite{dissertation} and \cite{compositionpaper} ultradifferentiable classes defined in terms of weight matrices $\mathcal{M}$ have been introduced and it has been shown that, by using the weight matrix $\Omega$ {\itshape associated with a given weight function} $\omega$, in this general framework one is able to treat both classical settings in a uniform and convenient way but also more classes.\vspace{6pt}

Since to each weight function $\omega$ we can associate a weight matrix $\Omega$ and since it is known that to each sequence $M$ one can associate a weight function $\omega_M$ (see \cite[Chapitre I]{mandelbrojtbook}), in \cite[Sect. 9.3]{dissertation} and in \cite{testfunctioncharacterization} the following iterated process has been studied: When $\mathcal{I}=\RR_{>0}$ is denoting the index set, then by starting with an abstractly given matrix $\mathcal{M}:=\{M^x: x\in\mathcal{I}\}$ with some regularity properties we immediately get the weight function matrix $\omega_{\mathcal{M}}:=\{\omega_{M^x}: x\in\mathcal{I}\}$. For each associated weight function $\omega_{M^x}$ one considers the associated weight matrix ${\Omega_{M^x}}$ and then proceed by iteration. This construction leads to multi-index weight matrices (and corresponding classes of ultradifferentiable functions) and the main questions in \cite[Sect. 9.3]{dissertation} and \cite{testfunctioncharacterization} have been the following:

\begin{itemize}
\item[$(*)$] Study the effects of growth properties assumed for $\mathcal{M}$ on this construction.

\item[$(**)$] Describe the class $\mathcal{E}_{[\mathcal{M}]}$ alternatively (as a locally convex vector space) by $\mathcal{E}_{[\omega_{\mathcal{M}}]}$. This enables the possibility to apply techniques from the classical (single) weight function setting to the space $\mathcal{E}_{[\mathcal{M}]}$.
\end{itemize}

Concerning $(**)$ we can see the following ''dual problem'':

\begin{itemize}
\item[$(**)'$] Starting with an abstractly given {\itshape weight function matrix} $\mathcal{M}_{\mathcal{W}}:=\{\omega^x: x\in\mathcal{I}\}$, can we describe the class $\mathcal{E}_{[\mathcal{M}_{\mathcal{W}}]}$ alternatively (as a locally convex vector space) by $\mathcal{E}_{[\mathcal{N}]}$ for some weight matrix $\mathcal{N}=\{N^x: x\in\mathcal{I}\}$ ?
\end{itemize}

For treating $(*)$ and $(**)$, in \cite[Sect. 9.3]{dissertation} and \cite{testfunctioncharacterization} the following properties for $\mathcal{M}$ have become relevant, see also \cite[Sect. 4.1]{compositionpaper} and \cite[Sect. 7.2]{dissertation}):\par\vskip.3cm

\hypertarget{R-mg}{$(\mathcal{M}_{\{\text{mg}\}})$} \hskip1cm $\forall\;x\in\mathcal{I}\;\exists\;C>0\;\exists\;y\in\mathcal{I}\;\forall\;j,k\in\NN:\;\;\; M^x_{j+k}\le C^{j+k} M^y_j M^y_k$,\par\vskip.3cm

\hypertarget{B-mg}{$(\mathcal{M}_{(\text{mg})})$} \hskip1cm $\forall\;x\in\mathcal{I}\;\exists\;C>0\;\exists\;y\in\mathcal{I}\;\forall\;j,k\in\NN:\;\;\; M^y_{j+k}\le C^{j+k} M^x_j M^x_k$,\par\vskip.3cm
\hypertarget{R-L}{$(\mathcal{M}_{\{\text{L}\}})$} \hskip1cm $\forall\;C>0\;\forall\;x\in\mathcal{I}\;\exists\;D>0\;\exists\;y\in\mathcal{I}\;\forall\;j\in\NN:\;\;\;C^j M^x_j\le D M^y_j$,\vskip.3cm
\hypertarget{B-L}{$(\mathcal{M}_{(\text{L})})$} \hskip1cm $\forall\;C>0\;\forall\;x\in\mathcal{I}\;\exists\;D>0\;\exists\;y\in\mathcal{I}\;\forall\;j\in\NN:\;\;\;C^j M^y_j\le D M^x_j$.\vskip.3cm

All these conditions are automatically valid for $\Omega$, when $\omega$ is satisfying standard properties. \hyperlink{R-mg}{$(\mathcal{M}_{\{\on{mg}\}})$} and \hyperlink{B-mg}{$(\mathcal{M}_{(\on{mg})})$} are the natural generalizations of condition {\itshape moderate growth} arising frequently in the weight sequence setting. In \cite[Sect. 9.3]{dissertation} and \cite{testfunctioncharacterization} only the sufficiency of \hyperlink{R-L}{$(\mathcal{M}_{\{\text{L}\}})$} and \hyperlink{B-L}{$(\mathcal{M}_{(\text{L})})$} has been applied for the study of question $(**)$ and these conditions seemed to be too strong, see also the discussion in $(ii)$ in Remark \ref{mixedmgandLremark}. Note that $(\mathcal{M}_{[\text{L}]})$ is indispensable to have $\mathcal{E}_{[\Omega]}=\mathcal{E}_{[\omega]}$, see \eqref{newexpabsorb}.\vspace{6pt}

The main aim of this article is to give a complete solution to $(**)$ and to $(**)'$ expressed in terms of growth properties for $\mathcal{M}$. In order to do so, instead of \hyperlink{R-L}{$(\mathcal{M}_{\{\text{L}\}})$} we consider the weaker requirement
\begin{equation}\label{mixedom1roum}
\forall\;x\in\mathcal{I}\;\exists\;y\in\mathcal{I}:\;\;\limsup_{t\rightarrow+\infty}\frac{\omega_{M^y}(2t)}{\omega_{M^x}(t)}<+\infty,
\end{equation}
and similarly
\begin{equation}\label{mixedom1beur}
\forall\;x\in\mathcal{I}\;\exists\;y\in\mathcal{I}:\;\;\limsup_{t\rightarrow+\infty}\frac{\omega_{M^x}(2t)}{\omega_{M^y}(t)}<+\infty,
\end{equation}
instead of \hyperlink{B-L}{$(\mathcal{M}_{(\text{L})})$}. \eqref{mixedom1roum} and \eqref{mixedom1beur} are the mixed versions of the standard assumption $\omega(2t)=O(\omega(t))$ as $t\rightarrow+\infty$ on weight functions (denoted by \hyperlink{om1}{$(\omega_1)$} in this work). These types of mixed conditions, for either weight sequences or weight functions, have frequently appeared in the literature in the study of mixed extension results in the ultradifferentiable or ultraholomorphic settings, see \cite{ChaumatChollet94}, \cite{surjectivity}, \cite{mixedramisurj}, \cite{mixedsectorialextensions}, \cite{mixedsectorialextensionsbeur}.\vspace{6pt}

In \cite{index} the main goal has been to give connections between standard/frequently used assumptions for weight sequences and weight functions in the ultradifferentiable and ultraholomorphic framework and the concept of O-regular variation. More precisely, this has been done for {\itshape moderate growth} for weight sequences and \hyperlink{om1}{$(\omega_1)$} for (associated) weight functions. Note that analogous growth properties are also showing up when dealing with different areas of weighted spaces in Functional Analysis, see the introduction in \cite{index}. We have been able to characterize these conditions in terms of growth properties, positivity and finiteness of weight indices.

Inspired by these known results and the mixed conditions mentioned before, for weight sequences and (their associated) weight functions we introduce in this paper new {\itshape mixed growth indices} which are related to mixed O-regular variation. Moreover, we study and compare these indices under the construction $M\mapsto\omega_M\mapsto\Omega_M$.\vspace{6pt}

Summarizing, it turns out that we can give an answer to problem $(**)$ by assuming \hyperlink{R-mg}{$(\mathcal{M}_{\{\on{mg}\}})$} and \eqref{mixedom1roum} resp. \hyperlink{B-mg}{$(\mathcal{M}_{(\on{mg})})$} and \eqref{mixedom1beur}, and problem $(**)'$ can by reduced to $(**)$. Alternatively, we have the possibility to express problem $(**)$ (and hence $(**)'$) purely in terms of these new concept of mixed indices. In particular, when considering $\mathcal{M}=\{M\}$ or $\mathcal{M}_{\mathcal{W}}=\{\omega\}$, then we get the results from \cite{BonetMeiseMelikhov07} expressed in the language of O-regular variation. Finally, by using these indices one is able to understand the difference between conditions \hyperlink{R-L}{$(\mathcal{M}_{\{\text{L}\}})$} and \eqref{mixedom1roum} resp. between \hyperlink{B-L}{$(\mathcal{M}_{(\text{L})})$} and \eqref{mixedom1beur} in a more quantitative way.\vspace{6pt}

The paper is structured as follows: In Section \ref{weightsandcond} we gather all relevant information on weight sequences, functions and matrices and we recall the multi-index construction. In Section \ref{mixedomega1section}, first we give an exhaustive characterization of \hyperlink{R-L}{$(\mathcal{M}_{\{\text{L}\}})$} and \hyperlink{B-L}{$(\mathcal{M}_{(\text{L})})$} resp. of \eqref{mixedom1roum} and \eqref{mixedom1beur} by means of growth conditions expressed in terms of the weight structures $\mathcal{M}$, $\omega_{\mathcal{M}}$ and of the multi-index construction, see Lemma \ref{Lrewriting} and Theorem \ref{omega1theorem}. It turns out that for the weight matrix $\Omega$ there is no difference between these requirements of the particular, Roumieu or Beurling, type, see Corollary \ref{Lrewritingcor}. Based on these characterizations, in Section \ref{mixedomega1indexsection} we introduce mixed growth indices for weight sequences and (associated) weight functions, see Proposition \ref{Lomega1comparison} and Theorem \ref{omega1importequivalence}.\vspace{6pt}

In Section \ref{mixedmgsection} the analogous procedure is done for the mixed moderate growth conditions \hyperlink{R-mg}{$(\mathcal{M}_{\{\on{mg}\}})$} and \hyperlink{B-mg}{$(\mathcal{M}_{(\on{mg})})$}, see Proposition \ref{mgomega6comparison} and Theorem \ref{mixedmgcharact}.\vspace{6pt}

In the final Section \ref{consequencesforultra} we apply the developed theory and prove the main results, Theorem \ref{Thm32testfuncnew} which is solving $(**)$, and Theorem \ref{generalweightmatrixthm} answering problem $(**)'$. The special (and classical) cases $\mathcal{M}=\{M\}$ and $\mathcal{M}_{\mathcal{W}}=\{\omega\}$ are discussed as well.

\section{Weights and conditions}\label{weightsandcond}
\subsection{General notation}
We write $\NN:=\{0,1,2,\dots\}$ and $\NN_{>0}:=\{1,2,3,\dots\}$.  For any real $x\ge 0$ we denote by $\lfloor x\rfloor$ the largest integer $j$ such that $j\le x$. With the symbol $\mathcal{E}$ we denote the class of all smooth functions.

\subsection{Weight sequences}\label{weightsequences}
Given a sequence $M=(M_j)_j\in\RR_{>0}^{\NN}$ we also use $\mu_j:=\frac{M_j}{M_{j-1}}$, $\mu_0:=1$. $M$ is called {\itshape normalized} if $1=M_0\le M_1$ holds true which can always be assumed without loss of generality.


$M$ is called {\itshape log-convex} if
$$\forall\;j\in\NN_{>0}:\;M_j^2\le M_{j-1} M_{j+1},$$
equivalently if $(\mu_j)_j$ is nondecreasing. If $M$ is log-convex and normalized, then both $j\mapsto M_j$ and $j\mapsto(M_j)^{1/j}$ are nondecreasing and $(M_j)^{1/j}\le\mu_j$ for all $j\in\NN_{>0}$.

For our purposes it is convenient to consider the following set of sequences
$$\hypertarget{LCset}{\mathcal{LC}}:=\{M\in\RR_{>0}^{\NN}:\;M\;\text{is normalized, log-convex},\;\lim_{j\rightarrow+\infty}(M_j)^{1/j}=+\infty\}.$$
We see that $M\in\hyperlink{LCset}{\mathcal{LC}}$ if and only if $(\mu_j)_j$ is nondecreasing and $\lim_{j\rightarrow+\infty}\mu_j=+\infty$ (e.g. see \cite[p. 104]{compositionpaper}). Moreover, there is a one-to-one correspondence between $M$ and $\mu=(\mu_j)_j$ by taking $M_j:=\prod_{i=0}^j\mu_i$.\vspace{6pt}

$M$ has the condition {\itshape moderate growth}, denoted by \hypertarget{mg}{$(\text{mg})$}, if
$$\exists\;C\ge 1\;\forall\;j,k\in\NN:\;M_{j+k}\le C^{j+k} M_j M_k.$$
In \cite{Komatsu73} it is denoted by $(M.2)$ and called {\itshape stability under ultradifferential operators.} 


Let $M,N\in\RR_{>0}^{\NN}$ be given, we write $M\hypertarget{preceq}{\preceq}N$ if $\sup_{j\in\NN_{>0}}\left(\frac{M_j}{N_j}\right)^{1/j}<+\infty$. We call $M$ and $N$ {\itshape equivalent}, denoted by $M\hypertarget{approx}{\approx}N$, if $M\hyperlink{preceq}{\preceq}N$ and $N\hyperlink{preceq}{\preceq}M$. Finally, write $M\le N$ if $M_j\le N_j$ for all $j\in\NN$. Note: If $M,N\in\RR_{>0}^{\NN}$ both are normalized, then $M\hyperlink{preceq}{\preceq}N$ does precisely mean $M_j\le C^jN_j$ for some $C\ge 1$ and all $j\in\NN$.

\subsection{Associated weight function}
Let $M\in\RR_{>0}^{\NN}$ (with $M_0=1$), then the {\itshape associated function} $\omega_M: \RR_{\ge 0}\rightarrow\RR\cup\{+\infty\}$ is defined by
\begin{equation*}\label{assofunc}
\omega_M(t):=\sup_{j\in\NN}\log\left(\frac{t^j}{M_j}\right)\;\;\;\text{for}\;t>0,\hspace{30pt}\omega_M(0):=0.
\end{equation*}
For an abstract introduction of the associated function we refer to \cite[Chapitre I]{mandelbrojtbook}, see also \cite[Definition 3.1]{Komatsu73}. If $\liminf_{j\rightarrow+\infty}(M_j)^{1/j}>0$, then $\omega_M(t)=0$ for sufficiently small $t$, since $\log\left(\frac{t^j}{M_j}\right)<0\Leftrightarrow t<(M_j)^{1/j}$ holds for all $j\in\NN_{>0}$. Moreover under this assumption $t\mapsto\omega_M(t)$ is a continuous nondecreasing function, which is convex in the variable $\log(t)$ and tends faster to infinity than any $\log(t^j)$, $j\ge 1$, as $t\rightarrow+\infty$. $\lim_{j\rightarrow+\infty}(M_j)^{1/j}=+\infty$ implies that $\omega_M(t)<+\infty$ for each finite $t$ which shall be considered as a basic assumption for defining $\omega_M$.\vspace{6pt}

If $M\le N$, then clearly $\omega_N\le\omega_M$ follows. If $M\in\hyperlink{LCset}{\mathcal{LC}}$, then we can compute $M$ by involving $\omega_M$ as follows, see \cite[Chapitre I, 1.4, 1.8]{mandelbrojtbook} (and also \cite[Prop. 3.2]{Komatsu73}):
\begin{equation}\label{Prop32Komatsu}
M_j=\sup_{t\ge 0}\frac{t^j}{\exp(\omega_{M}(t))},\;\;\;j\in\NN.
\end{equation}

\subsection{Weight functions}\label{weightfunctions}
A function $\omega:[0,+\infty)\rightarrow[0,+\infty)$ is called a {\itshape weight function} (in the terminology of \cite[Section 2.1]{index} and \cite[Section 2.2]{sectorialextensions}), if it is continuous, nondecreasing, $\omega(0)=0$ and $\lim_{t\rightarrow+\infty}\omega(t)=+\infty$. If $\omega$ satisfies in addition $\omega(t)=0$ for all $t\in[0,1]$, then we call $\omega$ a {\itshape normalized weight function}. For convenience we will write that $\omega$ has $\hypertarget{om0}{(\omega_0)}$ if it is a normalized weight.\vspace{6pt}


Moreover we consider the following conditions, this list of properties has already been used in ~\cite{dissertation}.

\begin{itemize}
\item[\hypertarget{om1}{$(\omega_1)}$] $\omega(2t)=O(\omega(t))$ as $t\rightarrow+\infty$, i.e. $\exists\;L\ge 1\;\forall\;t\ge 0:\;\;\;\omega(2t)\le L(\omega(t)+1)$.


\item[\hypertarget{om3}{$(\omega_3)$}] $\log(t)=o(\omega(t))$ as $t\rightarrow+\infty$ ($\Leftrightarrow\lim_{t\rightarrow+\infty}\frac{t}{\varphi_{\omega}(t)}=0$, $\varphi_{\omega}$ being the function defined next).

\item[\hypertarget{om4}{$(\omega_4)$}] $\varphi_{\omega}:t\mapsto\omega(e^t)$ is a convex function on $\RR$.


\item[\hypertarget{om6}{$(\omega_6)$}] $\exists\;H\ge 1\;\forall\;t\ge 0:\;2\omega(t)\le\omega(H t)+H$.



\end{itemize}


For convenience we define the sets
$$\hypertarget{omset0}{\mathcal{W}_0}:=\{\omega:[0,\infty)\rightarrow[0,\infty): \omega\;\text{has}\;\hyperlink{om0}{(\omega_0)},\hyperlink{om3}{(\omega_3)},\hyperlink{om4}{(\omega_4)}\},\hspace{20pt}\hypertarget{omset1}{\mathcal{W}}:=\{\omega\in\mathcal{W}_0: \omega\;\text{has}\;\hyperlink{om1}{(\omega_1)}\}.$$
For any $\omega\in\hyperlink{omset0}{\mathcal{W}_0}$ we define the {\itshape Legendre-Fenchel-Young-conjugate} of $\varphi_{\omega}$ by
\begin{equation}\label{legendreconjugate}
\varphi^{*}_{\omega}(x):=\sup\{x y-\varphi_{\omega}(y): y\ge 0\},\;\;\;x\ge 0,
\end{equation}
with the following properties, e.g. see \cite[Remark 1.3, Lemma 1.5]{BraunMeiseTaylor90}: It is convex and nondecreasing, $\varphi^{*}_{\omega}(0)=0$, $\varphi^{**}_{\omega}=\varphi_{\omega}$, $\lim_{x\rightarrow+\infty}\frac{x}{\varphi^{*}_{\omega}(x)}=0$ and finally $x\mapsto\frac{\varphi_{\omega}(x)}{x}$ and $x\mapsto\frac{\varphi^{*}_{\omega}(x)}{x}$ are nondecreasing on $[0,+\infty)$. Note that by normalization we can extend the supremum in \eqref{legendreconjugate} from $y\ge 0$ to $y\in\RR$ without changing the value of $\varphi^{*}_{\omega}(x)$ for given $x\ge 0$.

Let $\sigma,\tau$ be weight functions, we write $\sigma\hypertarget{ompreceq}{\preceq}\tau$ if $\tau(t)=O(\sigma(t))\;\text{as}\;t\rightarrow+\infty$
and call them equivalent, denoted by $\sigma\hypertarget{sim}{\sim}\tau$, if
$\sigma\hyperlink{ompreceq}{\preceq}\tau$ and $\tau\hyperlink{ompreceq}{\preceq}\sigma$.\vspace{6pt}

We recall the following known result, e.g. see \cite[Lemma 2.8]{testfunctioncharacterization} resp. \cite[Lemma 2.4 $(i)$]{sectorialextensions} and the references mentioned in the proofs there.

\begin{lemma}\label{assoweightomega0}
Let $M\in\hyperlink{LCset}{\mathcal{LC}}$, then $\omega_M\in\hyperlink{omset0}{\mathcal{W}_0}$ holds true.
\end{lemma}

\subsection{Weight matrices}\label{classesweightmatrices}
For the following definitions and conditions see also \cite[Section 4]{compositionpaper}.

Let $\mathcal{I}=\RR_{>0}$ denote the index set (equipped with the natural order), a {\itshape weight matrix} $\mathcal{M}$ associated with $\mathcal{I}$ is a (one parameter) family of weight sequences $\mathcal{M}:=\{M^x\in\RR_{>0}^{\NN}: x\in\mathcal{I}\}$, such that
$$\forall\;x\in\mathcal{I}:\;M^x\;\text{is normalized, nondecreasing},\;M^{x}\le M^{y}\;\text{for}\;x\le y.$$
We call a weight matrix $\mathcal{M}$ {\itshape standard log-convex,} denoted by \hypertarget{Msc}{$(\mathcal{M}_{\on{sc}})$}, if
$$\forall\;x\in\mathcal{I}:\;M^x\in\hyperlink{LCset}{\mathcal{LC}}.$$

A matrix is called {\itshape constant} if $M^x\hyperlink{approx}{\approx}M^y$ for all $x,y\in\mathcal{I}$.\vspace{6pt}

On the set of all weight matrices we consider the following relations, see \cite[p. 111]{compositionpaper}: We write $\mathcal{M}\{\preceq\}\mathcal{N}$, if
$$\forall\;x\in\mathcal{I}\;\exists\;y\in\mathcal{I}:\;\;\;M^x\hyperlink{preceq}{\preceq}N^y,$$
and
$\mathcal{M}(\preceq)\mathcal{N}$, if
$$\forall\;x\in\mathcal{I}\;\exists\;y\in\mathcal{I}:\;\;\;M^y\hyperlink{preceq}{\preceq}N^x.$$
We write $\mathcal{M}\{\approx\}\mathcal{N}$, if $\mathcal{M}\{\preceq\}\mathcal{N}$ and $\mathcal{N}\{\preceq\}\mathcal{M}$ and $\mathcal{M}(\approx)\mathcal{N}$, if $\mathcal{M}(\preceq)\mathcal{N}$ and $\mathcal{N}(\preceq)\mathcal{M}$. Finally we call $\mathcal{M}$ and $\mathcal{N}$ {\itshape equivalent}, if both relations $(\approx)$ and $\{\approx\}$ hold true.

\begin{remark}\label{mixedmgandLremark}
Concerning the conditions \hyperlink{R-mg}{$(\mathcal{M}_{\{\on{mg}\}})$}, \hyperlink{B-mg}{$(\mathcal{M}_{(\on{mg})})$} and \hyperlink{R-L}{$(\mathcal{M}_{\{\on{L}\}})$}, \hyperlink{B-L}{$(\mathcal{M}_{(\on{L})})$} we summarize (see also the discussion in Remark \ref{equivremark}):

\begin{itemize}
\item[$(i)$] If $\mathcal{M}=\{M\}$ resp. more generally if $\mathcal{M}=\{M^x: x\in\mathcal{I}\}$ is constant, then $\mathcal{M}$ does have \hyperlink{R-mg}{$(\mathcal{M}_{\{\on{mg}\}})$} and/or \hyperlink{B-mg}{$(\mathcal{M}_{(\on{mg})})$} if and only if $M$ does have \hyperlink{mg}{$(\on{mg})$} resp. some/each $M^x$ does have \hyperlink{mg}{$(\on{mg})$}.

    Moreover, it is clear that $(\mathcal{M}_{[\on{mg}]})$ is stable under relation $[\approx]$.

\item[$(ii)$] It is immediate to see that \hyperlink{R-L}{$(\mathcal{M}_{\{\on{L}\}})$} and/or \hyperlink{B-L}{$(\mathcal{M}_{(\on{L})})$} can never be valid for any $\mathcal{M}=\{M\}$. However, it is possible that both conditions are fulfilled for a constant matrix. Moreover, in both conditions we are interested in $C>1$ since for $0<C\le 1$ they are clearly valid with $x=y$ and $D=1$.

    In general it is not clear that $(\mathcal{M}_{[\on{L}]})$ is preserved under $[\approx]$.
\end{itemize}
\end{remark}

\subsection{Weight matrices obtained by weight functions and multi-index weight matrices}\label{weightmatrixfromfunction}
We summarize some facts which are shown in \cite[Section 5]{compositionpaper} and are needed in this work. All properties listed below are valid for $\omega\in\hyperlink{omset0}{\mathcal{W}_0}$, except \eqref{newexpabsorb} for which \hyperlink{om1}{$(\omega_1)$} is necessary.

\begin{itemize}
\item[$(i)$] The idea was that to each $\omega\in\hyperlink{omset0}{\mathcal{W}_0}$ we can associate a standard log-convex weight matrix $\Omega:=\{W^l=(W^l_j)_{j\in\NN}: l>0\}$ by\vspace{6pt}

    \centerline{$W^l_j:=\exp\left(\frac{1}{l}\varphi^{*}_{\omega}(lj)\right)$.}\vspace{6pt}


\item[$(ii)$] $\Omega$ satisfies
    \begin{equation}\label{newmoderategrowth}
    \forall\;l>0\;\forall\;j,k\in\NN:\;\;\;W^l_{j+k}\le W^{2l}_jW^{2l}_k,
    \end{equation}
    so both \hyperlink{R-mg}{$(\mathcal{M}_{\{\on{mg}\}})$} and \hyperlink{B-mg}{$(\mathcal{M}_{(\on{mg})})$} are satisfied.

\item[$(iii)$] In case $\omega$ has in addition \hyperlink{om1}{$(\omega_1)$}, then $\Omega$ has also both \hyperlink{R-L}{$(\mathcal{M}_{\{\on{L}\}})$} and \hyperlink{B-L}{$(\mathcal{M}_{(\on{L})})$}, more precisely
     \begin{equation}\label{newexpabsorb}
     \forall\;h\ge 1\;\exists\;A\ge 1\;\forall\;l>0\;\exists\;D\ge 1\;\forall\;j\in\NN:\;\;\;h^jW^l_j\le D W^{Al}_j.
     \end{equation}

\item[$(iv)$] Equivalent weight functions yield equivalent associated weight matrices.


\item[$(v)$] \hyperlink{om6}{$(\omega_6)$} holds if and only if some/each $W^l$ satisfies \hyperlink{mg}{$(\on{mg})$} if and only if $W^l\hyperlink{approx}{\approx}W^n$ for each $l,n>0$. Consequently \hyperlink{om6}{$(\omega_6)$} is characterizing the situation when $\Omega$ is constant.

\item[$(vi)$] $\omega\hyperlink{sim}{\sim}\omega_{W^l}$ for each $l>0$.
\end{itemize}

In particular, given  $M\in\hyperlink{LCset}{\mathcal{LC}}$ we denote by $\hypertarget{OmomM}{\Omega_M}$ the weight matrix associated with $\omega_M$.

In \cite[Section 2.7]{testfunctioncharacterization} the new {\itshape weight function matrix} $\omega_{\mathcal{M}}:=\{\omega_{M^x}: x\in\mathcal{I}\}$ has been considered and ultradifferentiable classes $\mathcal{E}_{\{\omega_{\mathcal{M}}\}}$ and $\mathcal{E}_{(\omega_{\mathcal{M}})}$ have been introduced.\vspace{6pt}

Let us recall now this approach. First let $\mathcal{M}:=\{M^x: x\in\mathcal{I}\}$ be \hyperlink{Msc}{$(\mathcal{M}_{\text{sc}})$}. By Lemma \ref{assoweightomega0} we get $\omega_{M^x}\in\hyperlink{omset0}{\mathcal{W}_0}$ for each $x\in\mathcal{I}$.

Based on this information, in \cite[Section 3.1]{testfunctioncharacterization} the following multi-index weight matrix construction has been introduced:

\begin{equation*}
M^x\mapsto\omega_{M^x}\mapsto M^{x;l_1}\mapsto\omega_{M^{x;l_1}}\mapsto\ M^{x;l_1,l_2}\mapsto\dots,
\end{equation*}
where for $x\in\mathcal{I}$, $l_j\in\RR_{>0}$, $j\in\NN_{>0}$, and $i\in\NN$ we put
\begin{equation*}
M^{x;l_1,\dots,l_{j+1}}_i:=\exp\left(\frac{1}{l_{j+1}}\varphi^{*}_{\omega_{M^{x;l_1,\dots,l_j}}}(l_{j+1}i)\right),\;\;\; M^{x;l_1}_i:=\exp\left(\frac{1}{l_1}\varphi^{*}_{\omega_{M^x}}(l_1i)\right),
\end{equation*}
respectively
$$\omega_{M^{x;l_1,\dots,l_j}}(t):=\sup_{p\in\NN}\log\left(\frac{t^p}{M_p^{x;l_1,\dots,l_j}}\right)\;\text{for}\;t>0,\;\hspace{20pt}\omega_{M^{x;l_1,\dots,l_j}}(0):=0.$$

In this notation, for any $x\in\mathcal{I}$ we get
\begin{equation}\label{Mxonestable}
M^{x;1}_j=M^x_j,\;\;\;j\in\NN,
\end{equation}
which follows by applying \eqref{Prop32Komatsu} (see also the proof of \cite[Thm. 6.4]{testfunctioncharacterization}):
     \begin{align*}
     M^{x;1}_j&:=\exp(\varphi^{*}_{\omega_{M}}(j))=\exp(\sup_{y\ge 0}\{jy-\omega_{M}(e^y)\})=\sup_{y\ge 0}\exp(jy-\omega_{M}(e^y))
     \\&
     =\sup_{y\ge 0}\frac{\exp(jy)}{\exp(\omega_{M}(e^y))}=\sup_{t\ge 1}\frac{t^j}{\exp(\omega_{M}(t))}=\sup_{t\ge 0}\frac{t^j}{\exp(\omega_{M}(t))}=M_j.
     \end{align*}
     Note that by normalization of $M$ we have $\omega_{M}(t)=0$ for $0\le t\le 1$, e.g. see the known integral representation formula for $\omega_M$, \cite[1.8. III]{mandelbrojtbook} and also \cite[$(3.11)$]{Komatsu73}. The similar statement is valid for the higher multi-index sequences.\vspace{6pt}

For a given \hyperlink{Msc}{$(\mathcal{M}_{\text{sc}})$} matrix $\mathcal{M}:=\{M^x: x\in\mathcal{I}\}$ we put $\mathcal{M}^{(2)}:=\{M^{x;l}: x\in\mathcal{I}, l>0\}$. An abstractly given matrix of weight functions denotes a family of weight functions $\{\omega^x: x\in\mathcal{I}\}$, such that $\omega^y\le\omega^x$ for any $x,y\in\mathcal{I}$ with $x\le y$.

\subsection{Ultradifferentiable classes}\label{classes}
Let $U\subseteq\RR^d$ be non-empty open. We write $K\subset\subset U$ if $K$ is compact and $K\subseteq U$. We introduce now the following spaces of ultradifferentiable functions classes.
First, for weight sequences we define the (local) classes of Roumieu type by
$$\mathcal{E}_{\{M\}}(U):=\{f\in\mathcal{E}(U):\;\forall\;K\subset\subset U\;\exists\;h>0:\;\|f\|_{M,K,h}<+\infty\},$$
and the classes of Beurling type by
$$\mathcal{E}_{(M)}(U):=\{f\in\mathcal{E}(U):\;\forall\;K\subset\subset U\;\forall\;h>0:\;\|f\|_{M,K,h}<+\infty\},$$
where we denote
\begin{equation*}
\|f\|_{M,K,h}:=\sup_{\alpha\in\NN^d,x\in K}\frac{|f^{(\alpha)}(x)|}{h^{|\alpha|} M_{|\alpha|}}.
\end{equation*}

For a compact set $K$ with sufficiently regular (smooth) boundary
$$\mathcal{E}_{M,h}(K):=\{f\in\mathcal{E}(K): \|f\|_{M,K,h}<+\infty\}$$
is a Banach space and so we have the following topological vector spaces
$$\mathcal{E}_{\{M\}}(K):=\underset{h>0}{\varinjlim}\;\mathcal{E}_{M,h}(K), \qquad \text{and} \qquad
\mathcal{E}_{\{M\}}(U)=\underset{K\subset\subset U}{\varprojlim}\;\underset{h>0}{\varinjlim}\;\mathcal{E}_{M,h}(K)=\underset{K\subset\subset U}{\varprojlim}\;\mathcal{E}_{\{M\}}(K).$$
Similarly, we get
$$\mathcal{E}_{(M)}(K):=\underset{h>0}{\varprojlim}\;\mathcal{E}_{M,h}(K),
\qquad \text{and} \qquad
\mathcal{E}_{(M)}(U)=\underset{K\subset\subset U}{\varprojlim}\;\underset{h>0}{\varprojlim}\;\mathcal{E}_{M,h}(K)=\underset{K\subset\subset U}{\varprojlim}\;\mathcal{E}_{(M)}(K).$$

For any $\omega\in\hyperlink{omset0}{\mathcal{W}_0}$ let the Roumieu type class be defined by
$$\mathcal{E}_{\{\omega\}}(U):=\{f\in\mathcal{E}(U):\;\forall\;K\subset\subset U\;\exists\;l>0:\;\|f\|_{\omega,K,l}<+\infty\}$$
and the space of Beurling type by
$$\mathcal{E}_{(\omega)}(U):=\{f\in\mathcal{E}(U):\;\forall\;K\subset\subset U\;\forall\;l>0:\;\|f\|_{\omega,K,l}<+\infty\},$$
where
\begin{equation*}
\|f\|_{\omega,K,l}:=\sup_{\alpha\in\NN^d,x\in K}\frac{|f^{(\alpha)}(x)|}{\exp(\frac{1}{l}\varphi^{*}_{\omega}(l|\alpha|))}.
\end{equation*}
For compact sets $K$ with sufficiently regular (smooth) boundary
$$\mathcal{E}_{\omega,l}(K):=\{f\in\mathcal{E}(K): \|f\|_{\omega,K,l}<+\infty\}$$
is a Banach space and we have the following topological vector spaces
$$\mathcal{E}_{\{\omega\}}(K):=\underset{l>0}{\varinjlim}\;\mathcal{E}_{\omega,l}(K), \qquad \text{and} \qquad
\mathcal{E}_{\{\omega\}}(U)=\underset{K\subset\subset U}{\varprojlim}\;\underset{l>0}{\varinjlim}\;\mathcal{E}_{\omega,l}(K)=\underset{K\subset\subset U}{\varprojlim}\;\mathcal{E}_{\{\omega\}}(K).$$
Similarly, we get
$$\mathcal{E}_{(\omega)}(K):=\underset{l>0}{\varprojlim}\;\mathcal{E}_{\omega,l}(K),\qquad \text{and} \qquad
\mathcal{E}_{(\omega)}(U)=\underset{K\subset\subset U}{\varprojlim}\;\underset{l>0}{\varprojlim}\;\mathcal{E}_{\omega,l}(K)=\underset{K\subset\subset U}{\varprojlim}\;\mathcal{E}_{(\omega)}(K).$$

Next, we consider classes defined by weight matrices of Roumieu type $\mathcal{E}_{\{\mathcal{M}\}}$ and of Beurling type $\mathcal{E}_{(\mathcal{M})}$ as follows, see also \cite[4.2]{compositionpaper}. For all $K\subset\subset U$ with sufficiently smooth boundary we put
\begin{equation*}
\mathcal{E}_{\{\mathcal{M}\}}(K):=\bigcup_{x\in\mathcal{I}}\mathcal{E}_{\{M^x\}}(K),\hspace{20pt}\mathcal{E}_{\{\mathcal{M}\}}(U):=\bigcap_{K\subset\subset U}\bigcup_{x\in\mathcal{I}}\mathcal{E}_{\{M^x\}}(K),
\end{equation*}
and
\begin{equation*}
\mathcal{E}_{(\mathcal{M})}(K):=\bigcap_{x\in\mathcal{I}}\mathcal{E}_{(M^x)}(K),\hspace{20pt}\mathcal{E}_{(\mathcal{M})}(U):=\bigcap_{x\in\mathcal{I}}\mathcal{E}_{(M^x)}(U).
\end{equation*}
For $K\subset\subset\RR^d$ one has the representation
$$\mathcal{E}_{\{\mathcal{M}\}}(K)=\underset{x\in\mathcal{I}}{\varinjlim}\;\underset{h>0}{\varinjlim}\;\mathcal{E}_{M^x,h}(K)$$
and so for $U\subseteq\RR^d$ non-empty open
\begin{equation*}
\mathcal{E}_{\{\mathcal{M}\}}(U)=\underset{K\subset\subset U}{\varprojlim}\;\underset{x\in\mathcal{I}}{\varinjlim}\;\underset{h>0}{\varinjlim}\;\mathcal{E}_{M^x,h}(K).
\end{equation*}
Similarly we get for the Beurling case
\begin{equation*}
\mathcal{E}_{(\mathcal{M})}(U)=\underset{K\subset\subset U}{\varprojlim}\;\underset{x\in\mathcal{I}}{\varprojlim}\;\underset{h>0}{\varprojlim}\;\mathcal{E}_{M^x,h}(K).
\end{equation*}
Finally, for the matrix $\omega_{\mathcal{M}}$ we put
$$\mathcal{E}_{\{\omega_{\mathcal{M}}\}}(U):=\{f\in\mathcal{E}(U):\;\forall\;K\subset\subset U\;\exists\;x\in\mathcal{I}\;\exists\;l>0:\;\|f\|_{\omega_{M^x},K,l}<+\infty\}$$
and
$$\mathcal{E}_{(\omega_{\mathcal{M}})}(U):=\{f\in\mathcal{E}(U):\;\forall\;K\subset\subset U\;\forall\;x\in\mathcal{I}\;\forall\;l>0:\;\|f\|_{\omega_{M^x},K,l}<+\infty\}.$$
Thus we obtain the topological vector spaces representations
\begin{equation*}
\mathcal{E}_{\{\omega_{\mathcal{M}}\}}(U)=\underset{K\subset\subset U}{\varprojlim}\;\underset{x\in\mathcal{I},l>0}{\varinjlim}\mathcal{E}_{\omega_{M^x},l}(K)
\end{equation*}
and
\begin{equation*}
\mathcal{E}_{(\omega_{\mathcal{M}})}(U)=\underset{K\subset\subset U}{\varprojlim}\;\underset{x\in\mathcal{I},l>0}{\varprojlim}\mathcal{E}_{\omega_{M^x},l}(K).
\end{equation*}
Since all arising limits are equal to countable ones, in each setting the Beurling type class is a Fr\'{e}chet space.

Ultradifferentiable classes of multi-index matrices can be defined in a similar way; for our purposes we will only need $\mathcal{M}^{(2)}$, see Section \ref{consequencesforultra}.

\section{The mixed $(\omega_1)$ conditions}\label{mixedomega1section}
Concerning \hyperlink{R-L}{$(\mathcal{M}_{\{\text{L}\}})$} and \hyperlink{B-L}{$(\mathcal{M}_{(\text{L})})$} we start with the following characterization.

\begin{lemma}\label{Lrewriting}
Let $\mathcal{M}:=\{M^x: x\in\mathcal{I}\}$ be \hyperlink{Msc}{$(\mathcal{M}_{\on{sc}})$} and let $\omega_{\mathcal{M}}:=\{\omega_{M^x}: x\in\mathcal{I}\}$ be the according matrix of associated weight functions. Then we get:
\begin{itemize}
\item[$(I)$]  The following conditions are equivalent:

\begin{itemize}
\item[$(i)$] $\mathcal{M}$ does satisfy \hyperlink{R-L}{$(\mathcal{M}_{\{\on{L}\}})$}.

\item[$(ii)$] The matrix $\mathcal{M}$ does satisfy
\begin{equation}\label{newLroum}
\forall\;x\in\mathcal{I}\;\exists\;D>0\;\exists\;y\in\mathcal{I}\;\forall\;j\in\NN:\;\;\;2^j M^x_j\le D M^y_j.
\end{equation}

\item[$(iii)$] The matrix $\omega_{\mathcal{M}}$ does satisfy
\begin{equation*}
\forall\; x\in\mathcal{I}\;\exists\;y\in\mathcal{I}\;\exists\; D>0\;\forall\;t\ge 0:\;\;\; \omega_{M^y}(2t)\le \omega_{M^x}(t)+D.
\end{equation*}

\item[$(iv)$] The matrix $\omega_{\mathcal{M}}$ does satisfy
\begin{equation*}
\forall\; C>0\;\forall\; x\in\mathcal{I}\;\exists\;y\in\mathcal{I}\;\exists\; D>0\;\forall\;t\ge 0:\;\;\; \omega_{M^y}(Ct)\le \omega_{M^x}(t)+D.
\end{equation*}
\end{itemize}

\item[$(II)$] The following conditions are equivalent:
\begin{itemize}
\item[$(i)$] $\mathcal{M}$ does satisfy \hyperlink{B-L}{$(\mathcal{M}_{(\on{L})})$},

\item[$(ii)$] $\mathcal{M}$ does satisfy
\begin{equation*}
\forall\;x\in\mathcal{I}\;\exists\;D>0\;\exists\;y\in\mathcal{I}\;\forall\;k\in\NN:\;\;\;2^j M^y_j\le D M^x_j.
\end{equation*}

\item[$(iii)$] The matrix $\omega_{\mathcal{M}}$ does satisfy
\begin{equation}\label{newLroummodomega2}
\forall\; x\in\mathcal{I}\;\exists\;y\in\mathcal{I}\;\exists\; D>0\;\forall\;t\ge 0:\;\;\; \omega_{M^x}(2t)\le \omega_{M^y}(t)+D.
\end{equation}

\item[$(iv)$] The matrix $\omega_{\mathcal{M}}$ does satisfy
\begin{equation*}
\forall\; C>0\;\forall\; x\in\mathcal{I}\;\exists\;y\in\mathcal{I}\;\exists\; D>0\;\forall\;t\ge 0:\;\;\; \omega_{M^x}(Ct)\le \omega_{M^y}(t)+D.
\end{equation*}
\end{itemize}
\end{itemize}
\end{lemma}

\demo{Proof}
We will only treat the Roumieu case, the Beurling case is analogous.\vspace{6pt}

$(i)\Rightarrow(ii)$ It suffices to take $C=2$ in condition \hyperlink{R-L}{$(\mathcal{M}_{\{\on{L}\}})$}.\vspace{6pt}

$(ii)\Rightarrow(iii)$ For all $x\in\mathcal{I}$ we can find $D>0$ and $y\in\mathcal{I}$ such that $\frac{(2t)^j}{M^y_j}\le D\frac{t^j}{M^y_j}$ for all $j\in\NN$ and $t\ge 0$. This yields the conclusion by definition of associated weight functions.\vspace{6pt}


$(iii)\Rightarrow(iv)$ This is clear for $C\le 2$ (since each associated weight function is nondecreasing). If $C>2$, then we take $n\in\NN$ chosen minimal such that $C\le 2^n$ is valid and apply iteration: Given $x\in\mathcal{I}$, there exists $y_1\in\mathcal{I}$ and $D_1>0$ such that for all $t\ge 0$,
$\omega_{M^{y_1}}(2t)\le \omega_{M^x}(t)+D_1$. Recursively, for $j=2,\dots,n$ we find there exist $y_j\in\mathcal{I}$ and $D_j>0$ such that for all $t\ge 0$,
$\omega_{M^{y_j}}(2t)\le \omega_{M^{y_{j-1}}}(t)+D_j$.
Then, we easily get
$$
\omega_{M^{y_n}}(Ct)\le \omega_{M^{y_n}}(2^nt)\le\omega_{M^x}(t)+\sum_{j=1}^n D_j,
$$
as desired.\vspace{6pt}

$(iv)\Rightarrow(i)$ We apply \eqref{Prop32Komatsu} and get for all $j\in\NN$:
$$M^x_j=\sup_{t\ge 0}\frac{t^j}{\exp(\omega_{M^x}(t))}\le e^D\sup_{t\ge 0}\frac{t^j}{\exp(\omega_{M^y}(Ct))}=e^D\frac{1}{C^j}\sup_{s\ge 0}\frac{s^j}{\exp(\omega_{M^y}(s))}=e^D\frac{1}{C^j}M^y_j,$$
and so we are done.
\qed\enddemo

For abstractly given weight matrices $\mathcal{M}$, a connection to mixed ''\hyperlink{om1}{$(\omega_1)$}-conditions'' has been established in \cite[Prop. 3.12, Cor. 3.15]{testfunctioncharacterization}. The different equivalent conditions in the following result make intervene both the associated weight functions $\omega_{M^x}$ and the matrices associated with the latter. In particular, we generalize the main characterizing result \cite[Thm. 3.1]{subaddlike} from the weight sequence to the weight matrix setting.

\begin{theorem}\label{omega1theorem}
Let $\mathcal{M}=\{M^x: x\in\mathcal{I}\}$ be \hyperlink{Msc}{$(\mathcal{M}_{\on{sc}})$} and let $\omega_{\mathcal{M}}:=\{\omega_{M^x}: x\in\mathcal{I}\}$ be the according matrix of associated weight functions and $\mathcal{M}^{(2)}:=\{M^{x;l}: x\in\mathcal{I}, l>0\}$. Then we get:

\begin{itemize}
\item[$(I)$] The following conditions are equivalent:

\begin{itemize}

\item[$(i)$] The matrix $\omega_{\mathcal{M}}$ satisfies
\begin{equation}\label{R-L-consequequstrong}
\forall\;h>1\;\forall\;x\in\mathcal{I}\;\exists\;y\in\mathcal{I}:
\;\limsup_{t\rightarrow+\infty}\frac{\omega_{M^y}(ht)}{\omega_{M^x}(t)}<+\infty.
\end{equation}

\item[$(ii)$] The matrix $\omega_{\mathcal{M}}$ satisfies
\begin{equation}\label{R-L-consequ}
\forall\;x\in\mathcal{I}\;\exists\;y\in\mathcal{I}:\;\;
\limsup_{t\rightarrow+\infty}\frac{\omega_{M^y}(2t)}{\omega_{M^x}(t)}<+\infty.
\end{equation}

\item[$(iii)$] The matrix $\mathcal{M}$ does satisfy

\begin{equation}\label{omega1mixedcharactequmod}
\exists\;r>1\;\forall\;x\in\mathcal{I}\;\exists\;y\in\mathcal{I}\;\exists\;L\in\NN_{>0}:\;\;\;\liminf_{j\rightarrow+\infty}\frac{(M^y_{Lj})^{1/(Lj)}}{(M^x_j)^{1/j}}>r.
\end{equation}

\item[$(iv)$] The matrix $\mathcal{M}$ does satisfy

\begin{equation}\label{omega1mixedcharactequmodstrong}
\forall\;r>1\;\forall\;x\in\mathcal{I}\;\exists\;y\in\mathcal{I}\;\exists\;L\in\NN_{>0}:\;\;\;\liminf_{j\rightarrow+\infty}\frac{(M^y_{Lj})^{1/(Lj)}}{(M^x_j)^{1/j}}>r.
\end{equation}

\item[$(v)$] The matrix $\mathcal{M}^{(2)}$ does satisfy
\begin{equation*}
\forall\;x\in\mathcal{I}\;\forall\;C> 1\;\exists\;y\in\mathcal{I}\;\exists\;B>0\;\forall\;a>0\;\exists\;D>0\;\forall\;j\in\NN:\;\;\; C^jM_j^{x;a}\le D M_j^{y;Ba}.
\end{equation*}

\item[$(vi)$] The matrix $\mathcal{M}^{(2)}$ does satisfy
\begin{equation}\label{naturalitycor1roumnewnew}
\forall\;x\in\mathcal{I}\;\forall\;C> 1\;\forall\;a>0\;\exists\;y\in\mathcal{I}\;\exists\;b>0\;\exists\;D>0\;\forall\;j\in\NN:\;\;\;C^jM^{x;a}_j\le DM^{y;b}_j.
\end{equation}

\item[$(vii)$] The matrix $\mathcal{M}^{(2)}$ does satisfy
\begin{equation}\label{naturalitycor1new1}
   \forall\;r>1\;\forall\;x\in\mathcal{I}\;\exists\;y\in\mathcal{I}\;\exists\;b>0:\;\;\; \liminf_{j\in\NN_{>0}}\frac{\exp(\frac{1}{bj}\varphi^{*}_{\omega_{M^y}}(bj))}{\exp(\frac{1}{j}\varphi^{*}_{\omega_{M^x}}(j))}>r.
\end{equation}

\item[$(viii)$] The matrix $\mathcal{M}^{(2)}$ does satisfy
\begin{equation}\label{naturalitycor1new1real}
   \forall\;r>1\;\forall\;x\in\mathcal{I}\;\exists\;y\in\mathcal{I}\;\exists\;b>0:\;\;\; \liminf_{j\in\RR_{>0}}\frac{\exp(\frac{1}{bj}\varphi^{*}_{\omega_{M^y}}(bj))}{\exp(\frac{1}{j}\varphi^{*}_{\omega_{M^x}}(j))}>r.
\end{equation}
\end{itemize}

\item[$(II)$]  The following conditions are equivalent:

\begin{itemize}
\item[$(i)$] The matrix $\omega_{\mathcal{M}}$ satisfies
\begin{equation*}
\forall\;h>1\;\forall\;x\in\mathcal{I}\;\exists\;y\in\mathcal{I}:
\;\limsup_{t\rightarrow+\infty}\frac{\omega_{M^x}(ht)}{\omega_{M^y}(t)}<+\infty.
\end{equation*}

\item[$(ii)$] The matrix $\omega_{\mathcal{M}}$ satisfies
\begin{equation}\label{B-L-consequequ}
\forall\;x\in\mathcal{I}\;\exists\;y\in\mathcal{I}:\;\;
\limsup_{t\rightarrow+\infty}\frac{\omega_{M^x}(2t)}{\omega_{M^y}(t)}<+\infty.
\end{equation}

\item[$(iii)$] The matrix $\mathcal{M}$ does satisfy
\begin{equation}\label{omega1mixedcharactequbeurmod}
\exists\;r>1\;\forall\;x\in\mathcal{I}\;\exists\;y\in\mathcal{I}\;\exists\;L\in\NN_{>0}:\;\;\;\liminf_{j\rightarrow+\infty}\frac{(M^x_{Lj})^{1/(Lj)}}{(M^y_j)^{1/j}}>r.
\end{equation}

\item[$(iv)$] The matrix $\mathcal{M}$ does satisfy
\begin{equation}\label{omega1mixedcharactequbeurmodstrong}
\forall\;r>1\;\forall\;x\in\mathcal{I}\;\exists\;y\in\mathcal{I}\;\exists\;L\in\NN_{>0}:\;\;\;\liminf_{j\rightarrow+\infty}\frac{(M^x_{Lj})^{1/(Lj)}}{(M^y_j)^{1/j}}>r.
\end{equation}

\item[$(v)$] The matrix $\mathcal{M}^{(2)}$ does satisfy
\begin{equation*}
\forall\;x\in\mathcal{I}\;\forall\;C> 1\;\exists\;y\in\mathcal{I}\;\exists\;B>0\;\forall\;a>0\;\exists\;D>0\;\forall\;j\in\NN:\;\;\; C^jM_j^{y;a/B}\le D M_j^{x;a}.
\end{equation*}

\item[$(vi)$] The matrix $\mathcal{M}^{(2)}$ does satisfy
\begin{equation}\label{naturalitycor1roumnewnewbeur}
\forall\;x\in\mathcal{I}\;\forall\;C> 1\;\forall\;a>0\;\exists\;y\in\mathcal{I}\;\exists\;b>0\;\exists\;D>0\;\forall\;j\in\NN:\;\;\;C^jM^{y;b}_j\le DM^{x;a}_j.
\end{equation}

\item[$(vii)$] The matrix $\mathcal{M}^{(2)}$ does satisfy

\begin{equation*}
   \forall\;r>1\;\forall\;x\in\mathcal{I}\;\exists\;y\in\mathcal{I}\;\exists\;b>0:\;\;\; \liminf_{j\in\NN_{>0}}\frac{\exp(\frac{1}{j}\varphi^{*}_{\omega_{M^x}}(j))}{\exp(\frac{b}{j}\varphi^{*}_{\omega_{M^y}}(j/b))}>r.
\end{equation*}

\item[$(viii)$] The matrix $\mathcal{M}^{(2)}$ does satisfy
\begin{equation*}
   \forall\;r>1\;\forall\;x\in\mathcal{I}\;\exists\;y\in\mathcal{I}\;\exists\;b>0:\;\;\; \liminf_{j\in\RR_{>0}}\frac{\exp(\frac{1}{j}\varphi^{*}_{\omega_{M^x}}(j))}{\exp(\frac{b}{j}\varphi^{*}_{\omega_{M^y}}(j/b))}>r.
\end{equation*}
\end{itemize}
\end{itemize}
\end{theorem}

In particular, this result can be applied to any $\mathcal{M}=\{M^x: x\in\mathcal{I}\}$ being \hyperlink{Msc}{$(\mathcal{M}_{\on{sc}})$} and satisfying \hyperlink{R-L}{$(\mathcal{M}_{\{\on{L}\}})$} resp. \hyperlink{B-L}{$(\mathcal{M}_{(\on{L})})$}: In this case $(I)(i)\Leftrightarrow(iii)$ resp. $(II)(i)\Leftrightarrow(iii)$ in Lemma \ref{Lrewriting} yields that \eqref{R-L-consequ} resp. \eqref{B-L-consequequ}, i.e. the assertion $(ii)$ of each particular type, is valid.\vspace{6pt}

\demo{Proof}
Again, we treat the Roumieu case in detail, the Beurling setting is completely analogous.\vspace{6pt}

$(i)\Rightarrow(ii)$ This is clear.\vspace{6pt}

$(ii)\Rightarrow(iii)$ Let $x\in\mathcal{I}$ be given, so $\omega_{M^y}(2t)\le L\omega_{M^x}(t)+L$ for some index $y\in\mathcal{I}$, $L\in\NN_{>0}$ and all $t\ge 0$. Then, by applying \eqref{Prop32Komatsu}, we get for all $j\in\NN$:
\begin{align*}
M^y_{Lj}&=\sup_{t\ge 0}\frac{t^{Lj}}{\exp(\omega_{M^y}(t))}=\sup_{s\ge 0}\frac{(2s)^{Lj}}{\exp(\omega_{M^y}(2s))}\ge e^{-L}\sup_{s\ge 0}\frac{(2s)^{Lj}}{\exp(L\omega_{M^x}(s))}
\\&
=e^{-L}2^{Lj}\left(\sup_{s\ge 0}\frac{s^j}{\exp(\omega_{M^x}(s))}\right)^{L}=e^{-L}2^{Lj}(M^x_j)^{L},
\end{align*}
which proves $(iii)$.\vspace{6pt}

$(iii)\Rightarrow(iv)$ Let $r_0>1$ be the value given in $(iii)$. If $r\in(1,r_0]$, then nothing is to prove. If $r>r_0$, then we choose $n\in\NN$ minimal to have $r\le r_0^n$ and apply $n$ iterations. Namely, given $x\in\mathcal{I}$, there exists $y_1\in\mathcal{I}$ and $L_1\in\NN_{>0}$ such that
$$
\liminf_{j\rightarrow+\infty}\frac{(M^{y_1}_{L_1j})^{1/(L_1j)}}{(M^x_j)^{1/j}}>r_0.
$$
Recursively, for $i=2,\dots,n$ there exist $y_i\in\mathcal{I}$ and $L_i\in\NN_{>0}$ such that
$$
\liminf_{j\rightarrow+\infty}\frac{(M^{y_i}_{L_ij})^{1/(L_ij)}}{(M^{y_{i-1}}_j)^{1/j}}>r_0.
$$
The choice $y=y_n$ and $L=\prod_{i=1}^n L_i$ clearly fulfills the requirements.
\vspace{6pt}

$(iv)\Rightarrow(v)$ Recall that the matrix  $\{M^{x;a}: a>0\}$  associated with the weight $\omega_{M^x}$ is given by \begin{equation}\label{eq.Mxap}
M^{x;a}_j=\exp\left(\frac{1}{a}\varphi^{*}_{\omega_{M^x}}(a j)\right),\;\;\;a>0,\;\;\;j\in\NN,
\end{equation}
and $M^x\equiv M^{x;1}$ by \eqref{Mxonestable}.
Moreover, if $L\in\NN_{>0}$ we have for every $j\in\NN$,
\begin{equation*}
M^{x;L}_j=\exp\left(\frac{1}{L}\varphi^{*}_{\omega_{M^x}}(L j)\right)=(M^{x;1}_{Lj})^{1/L}=(M^x_{Lj})^{1/L}.
\end{equation*}
Let $x\in\mathcal{I}$ and $C>1$ be given. By the assumption (reasoning with $r=C$) we have
$$C^jM^{x;1}_j=C^jM^x_j\le A(M^y_{Lj})^{1/L}=AM^{y;L}_j$$
for some $y\in\mathcal{I}$, $L\in\NN_{>0}$, $A\ge 1$, and all $j\in\NN$.
Hence,
$$\log\left(\frac{(Ct)^j}{M^{y;L}_j}\right)\le \log\left(\frac{t^j}{M^{x;1}_j}\right)+\log(A)$$
for all $t>0$ and $j\in\NN$, and we obtain by definition $\omega_{M^{y;L}}(Ct)\le\omega_{M^{x;1}}(t)+\log(A)$ for all $t\ge 0$ (observe that $\omega_{M^{y;L}}(0)=\omega_{M^{x;1}}(0)=0$).
Recall that, as shown in \cite[Lemma 5.7]{compositionpaper} (see also \cite[Theorem 4.0.3, Lemma 5.1.3]{dissertation} and \cite[Lemma 2.5]{sectorialextensions} in a more precise way), we have for any $x\in\mathcal{I}$:
\begin{equation}\label{goodequivalence}
\forall\;a>0\;\exists\;D_a>0\;\forall\;t\ge 0:\;\;\;a\omega_{M^{x;a}}(t)\le\omega_{M^{x;1}}(t)=\omega_{M^x}(t)\le 2a\omega_{M^{x;a}}(t)+D_a.
\end{equation}
So we combine everything to get for all $t\ge 0$:
$$\omega_{M^y}(Ct)\le 2L\omega_{M^{y;L}}(Ct)+D_L\le 2L\omega_{M^{x;1}}(t)+2L\log(A)+D_L=2L\omega_{M^x}(t)+2L\log(A)+D_L.$$
From here, for every $s\in\RR$ we get
$$
\varphi_{\omega_{M^y}}(s+\log(C))=\omega_{M^y}(Ce^s)\le
2L\omega_{M^x}(e^s)+2L\log(A)+D_L=2L\varphi_{\omega_{M^x}}(s)+2L\log(A)+D_L,
$$
and consequently, for $a>0$ and every $j\in\NN$,
\begin{align*}
\varphi^*_{\omega_{M^x}}(aj)&=\sup_{t\ge 0}\{ajt-\varphi_{\omega_{M^x}}(t)\}\\
&\le \sup_{t\ge 0}\{aj(t+\log(C))-\frac{1}{2L}\varphi_{\omega_{M^y}}(t+\log(C))\}+\log(A)+\frac{D_L}{2L}-aj\log(C)\\
&\le \frac{1}{2L}\sup_{s\ge 0}\{2aLjs-\varphi_{\omega_{M^y}}(s)\}+ \log(A)+\frac{D_L}{2L}-aj\log(C)\\
&=\frac{1}{2L}\varphi^*_{\omega_{M^y}}(2aLj)+\log(A)+\frac{D_L}{2L}-aj\log(C).
\end{align*}
In conclusion, we use~\eqref{eq.Mxap} in order to obtain
$$
M^{x;a}_j\le A^{1/j}e^{D_L/(2aL)}\frac{1}{C^j}M^{y;2aL}_j,
$$
as desired.\vspace{6pt}

$(v)\Rightarrow(vi)$ It suffices to take, for any $a>0$, $b=Ba$, where $B$ is given in $(v)$.\vspace{6pt}

$(vi)\Rightarrow(vii)$ If we choose $a=1$ in \eqref{naturalitycor1roumnewnew} we get
$$
\forall\;x\in\mathcal{I}\;\forall\;C> 1\;\exists\;y\in\mathcal{I}\;\exists\;b>0\;\exists\;D>0\;\forall\;j\in\NN:\;\;\;C^jM^{x;1}_j\le DM^{y;b}_j,
$$
what leads to the conclusion in view of~\eqref{eq.Mxap}.\vspace{6pt}

$(vii)\Rightarrow(viii)$

Note that the mapping $j\mapsto\frac{1}{jb}\varphi^{*}_{\omega_{M^x}}(jb)$ is nondecreasing for any fixed $b>0$ and $x\in\mathcal{I}$. Let $j_1\in\RR_{>1}$ be given and take $j\in\NN_{\ge 2}$ with $j-1<j_1\le j$ and then
\begin{align*}
&\frac{\exp(\frac{1}{bj}\varphi^{*}_{\omega_{M^y}}(bj))}{\exp(\frac{1}{j}\varphi^{*}_{\omega_{M^x}}(j))}\le\frac{\exp(\frac{1}{bj}\varphi^{*}_{\omega_{M^y}}(bj))}{\exp(\frac{1}{j_1}\varphi^{*}_{\omega_{M^x}}(j_1))}\le\frac{\exp(\frac{1}{2bj_1}\varphi^{*}_{\omega_{M^y}}(2bj_1))}{\exp(\frac{1}{j_1}\varphi^{*}_{\omega_{M^x}}(j_1))},
\end{align*}
since $bj\le 2bj_1\Leftrightarrow j\le 2j_1$ is valid by $j\le 2(j-1)\Leftrightarrow 2\le j$. Thus we have verified \eqref{naturalitycor1new1real} with the same value $r>1$ and with the same choice $y\in\mathcal{I}$ for given index $x$ when taking $b':=2b$, $b>0$ denoting the parameter in \eqref{naturalitycor1new1}.\vspace{6pt}

$(viii)\Rightarrow(i)$ Given $h>1$ and $x\in\mathcal{I}$, we apply the hypothesis with some $r>h$ and deduce that
$h^jM^x_j\le AM^{y;b}_j$
for some $y\in\mathcal{I}$, $b>0$, $A\ge 1$, and all $j\in\NN$.
Hence,
$$\log\left(\frac{(ht)^j}{M^{y;b}_j}\right)\le \log\left(\frac{t^j}{M^{x}_j}\right)+\log(A)$$
for all $t>0$ and $j\in\NN$, and we obtain $\omega_{M^{y;b}}(ht)\le\omega_{M^{x}}(t)+\log(A)$ for all $t\ge 0$.
We conclude by using~\eqref{goodequivalence} again.
\qed\enddemo

We gather now the information for weight matrices associated with given Braun-Meise-Taylor weight functions $\omega$ and get the following characterization.

\begin{corollary}\label{Lrewritingcor}
Let $\omega\in\hyperlink{omset0}{\mathcal{W}_0}$ be given and let $\Omega=\{W^x: x>0\}$ be the associated weight matrix. Moreover let $\omega_{\Omega}:=\{\omega_{W^x}: x>0\}$ be the matrix of the associated weight functions. Then the following are equivalent:
\begin{itemize}
\item[$(i)$] $\omega$ satisfies \hyperlink{om1}{$(\omega_1)$}.

\item[$(ii)$] $\Omega$ satisfies \hyperlink{R-L}{$(\mathcal{M}_{\{\on{L}\}})$}.

\item[$(iii)$] $\omega_{\Omega}$, $\Omega$ and $\Omega^{(2)}$ satisfy any of the corresponding equivalent Roumieu-like conditions in $(I)$ in Theorem \ref{omega1theorem}.

\item[$(iv)$] $\Omega$ satisfies \hyperlink{B-L}{$(\mathcal{M}_{(\on{L})})$}.

\item[$(v)$] $\omega_{\Omega}$, $\Omega$ and $\Omega^{(2)}$ satisfy any of the corresponding equivalent Beurling-like conditions in $(II)$ in Theorem \ref{omega1theorem}.
\end{itemize}

So, the mixed \hyperlink{om1}{$(\omega_1)$} conditions of Roumieu and/or Beurling type are equivalent to \hyperlink{om1}{$(\omega_1)$} for $\omega$.
\end{corollary}

\demo{Proof}
$(i)\Leftrightarrow(iii)$, $(i)\Leftrightarrow(v)$ Both equivalences hold true by taking into account that  $\omega\hyperlink{sim}{\sim}\omega_{W^x}$ for all $x>0$, see \cite[Lemma 5.7]{compositionpaper}.\vspace{6pt}

$(i)\Rightarrow(ii),(iv)$ holds by \cite[Lemma 5.9 $(5.10)$]{compositionpaper}, see \eqref{newexpabsorb}.\vspace{6pt}

$(ii)\Rightarrow(iii)$, $(iv)\Rightarrow(v)$ Both implications follow by Lemma \ref{Lrewriting}.
\qed\enddemo

\begin{remark}\label{omega1mixedcharactrem}
We summarize now some facts concerning the arising $\liminf$-conditions in the previous results.

\begin{itemize}
\item[$(a)$] \eqref{newLroum} implies \eqref{omega1mixedcharactequmod} (with $L=1$), but the latter condition is weaker than the first one since $j\mapsto(M^x_j)^{1/j}$ is nondecreasing by log-convexity and normalization for each index $x$ fixed.

\item[$(b)$] This observation is consistent with the characterizations shown in Lemma \ref{Lrewriting} and Theorem \ref{omega1theorem}; more precisely compare the mixed \hyperlink{om1}{$(\omega_1)$} conditions \eqref{newLroummodomega2} and \eqref{R-L-consequ}. In the latter one on the right hand side we only require a $O$-growth restriction (and similarly for the Beurling setting).

\item[$(c)$] However, Corollary \ref{Lrewritingcor} yields that \eqref{newLroum} is equivalent to \eqref{omega1mixedcharactequmod} when the matrix $\mathcal{M}\equiv\Omega$ is associated with a given weight function $\omega\in\hyperlink{omset0}{\mathcal{W}_0}$.

\item[$(d)$] Concerning Corollary \ref{Lrewritingcor} we also make the following observation: When considering the matrix $\Omega$, then w.l.o.g. in $(I)(iii)$ and $(I)(iv)$  resp. in $(II)(iii)$ and $(II)(iv)$ in Theorem~\ref{omega1theorem} we can choose $y=x$ which can be seen directly as follows:

By definition of the associated weight matrix and the properties for $\varphi^{*}_{\omega}$ (see Section \ref{weightmatrixfromfunction}) we have
$$\forall\;y\ge x>0\;\forall\;L,L'\in\NN_{>0}, L'\ge\frac{yL}{x}\;\forall\;j\in\NN_{>0}:\;\;\;(W^y_{Lj})^{1/(Lj)}\le(W^{x}_{L'j})^{1/(L'j)}.$$
From this, the desired statement follows for the Roumieu type immediately.

In the Beurling case, when $0<y<x$, then we choose $L_1\in\NN_{>0}$ such that $L_1\ge\frac{x}{y}$ and so $(W^x_{j})^{1/j}\le(W^y_{L_1j})^{1/(L_1j)}$ holds true for any $j\in\NN_{>0}$. Thus we can estimate by
$$\liminf_{j\rightarrow\infty}\frac{(W^x_{Lj})^{1/(Lj)}}{(W^y_j)^{1/j}}\le\liminf_{j\rightarrow\infty}\frac{(W^x_{LL_1j})^{1/(LL_1j)}}{(W^y_{L_1j})^{1/(L_1j)}}\le\liminf_{j\rightarrow\infty}\frac{(W^x_{LL_1j})^{1/(LL_1j)}}{(W^x_{j})^{1/j}},$$
verifying the desired statement for the Beurling case as well.

\item[$(e)$] A different but more involved argument for the proof of $(d)$ is to combine Corollary \ref{Lrewritingcor} with $(vi)$ in Sect. \ref{weightmatrixfromfunction} and the main result \cite[Thm. 3.1]{subaddlike}: $\omega$ satisfies \hyperlink{om1}{$(\omega_1)$} if and only if some/each $\omega_{W^x}$ does, and this is equivalent to the fact that $(I)(iii)$ resp. $(II)(iii)$ in Theorem~\ref{omega1theorem} has to hold for some/each $x=y$.
\end{itemize}
\end{remark}

\subsection{Mixed growth indices based on mixed $(\omega_1)$ conditions}\label{mixedomega1indexsection}
We start this section with some observations. In \eqref{omega1mixedcharactequmod} and \eqref{omega1mixedcharactequbeurmod} in the arising $\liminf$ condition in order to make sense we have to assume $L\in\NN_{>0}$, whereas only $b>0$ is required in the conditions from $(vii)$ and $(viii)$ in Theorem \ref{omega1theorem}.

Moreover recall that, since $M^{x;a}\in\hyperlink{LCset}{\mathcal{LC}}$, by \eqref{Prop32Komatsu} we get
$$\forall\;x\in\mathcal{I}\;\forall\;a>0\;\forall\;j\in\NN:\;\;\;M^{x;a}_j=\sup_{t>0}\frac{t^j}{\exp(\omega_{M^{x;a}}(t))}=\exp(\varphi^{*}_{\omega_{M^{x;a}}}(j)),$$
and the last expression makes even sense for any real $j\ge 0$.\vspace{6pt}

Based on the characterizations shown in Theorem \ref{omega1theorem} we introduce now two mixed growth indices. Let $M,N\in\hyperlink{LCset}{\mathcal{LC}}$ with $M\le N$ and let $\omega,\sigma$ be weight functions with $\sigma\ge\omega$.

For $b>0$ we write $(M,\Omega_N)_{L,b}$, if
\begin{equation}\label{om1a}
\exists\;q>1:\;\liminf_{j\rightarrow+\infty}\frac{\exp(\frac{1}{jq}\varphi^{*}_{\omega_{N}}(jq))}{\exp(\frac{1}{j}\varphi^{*}_{\omega_{M}}(j))}>q^b,
\end{equation}
which shall be compared with \cite[Thm. 3.11 $(v)$]{index}.

Note that, if we write $\hyperlink{OmomM}{\Omega_N}=\{W^q=(W_j^q)_{j\in\NN}: q>0\}$ for the weight matrix associated with $\omega_N$, so in particular $W^1=N$ (see Section \ref{weightmatrixfromfunction}), then the previous condition reads
$$
\exists\;q>1:\;\liminf_{j\rightarrow+\infty}\left(\frac{W_j^q}{M_j}\right)^{1/j}>q^b,
$$
what explains the notation used. It is immediate to see that if the condition is satisfied for some given $b>0$, then also for all $0<b'<b$ (with the same choice for $q$).

Similarly, given $a>0$ we write $(\sigma,\omega)_{\omega_1,a}$ if
\begin{equation*}
\exists\;K>1:\;\limsup_{t\rightarrow+\infty}\frac{\omega(Kt)}{\sigma(t)}<K^{a},
\end{equation*}
which shall be compared with \cite[Thm. 2.11 $(iv)$]{index}. Again it is clear that if the condition is satisfied for some given $a>0$, then also for all $a'>a$ (with the same choice of $K$). According to these growth restrictions we put
\begin{equation}\label{om1a2}
\beta(M,\Omega_N):=\sup\{b>0: (M,\Omega_N)_{L,b}\},
\end{equation}
\begin{equation}\label{om1a3}
\alpha(\sigma,\omega):=\inf\{a>0: (\sigma,\omega)_{\omega_1,a}\}.
\end{equation}

If there does not exist any $b>0$, resp. $a>0$, such that $(M,\Omega_N)_{L,b}$, resp. $(\sigma,\omega)_{\omega_1,a}$, holds true, then we put $\beta(M,\Omega_N)=0$, resp. $\alpha(\sigma,\omega)=\infty$.

A first immediate consequence is the following:

\begin{lemma}\label{indexinfinity}
Let $M,N\in\hyperlink{LCset}{\mathcal{LC}}$ be given with $M\le N$ and satisfying
\begin{equation}\label{indexinfinityequ}
\exists\;C>1\;\exists\;D>0\;\forall\;j\in\NN:\;\;\;C^j M_j\le D N_j.
\end{equation}
Then $\alpha(\omega_M,\omega_N)=0$ holds true.
\end{lemma}

Note that \eqref{indexinfinityequ} cannot be valid for $M=N$ (see Remark \ref{mixedmgandLremark}) and in particular this result applies for any (non-constant) weight matrix $\mathcal{M}:=\{M^x: x\in\mathcal{I}\}$ being \hyperlink{Msc}{$(\mathcal{M}_{\on{sc}})$} and satisfying \hyperlink{R-L}{$(\mathcal{M}_{\{\on{L}\}})$} resp. \hyperlink{B-L}{$(\mathcal{M}_{(\on{L})})$}.

\demo{Proof}
The argument in $(ii)\Rightarrow(iii)$ in Lemma \ref{Lrewriting} may be repeated verbatim in order to yield that $\omega_{N}(Ct)\le\omega_{M}(t)+D_1$ for some $D_1\ge 1$ and all $t\ge 0$. Hence
$$\limsup_{t\rightarrow\infty}\frac{\omega_{N}(Ct)}{\omega_M(t)}\le 1<C^{a}$$
for any $a>0$, which gives the conclusion.
\qed\enddemo

\vspace{6pt}
The next main result shows that the values defined in \eqref{om1a2} and \eqref{om1a3} are closely related.

\begin{proposition}\label{Lomega1comparison}
Let $M,N\in\hyperlink{LCset}{\mathcal{LC}}$ be given with $M\le N$ and let $\omega_{M}, \omega_{N}$ be the corresponding associated weight functions, then we get
$$
\beta(M,\Omega_N)=\frac{1}{\alpha(\omega_{M},\omega_{N})}.
$$
\end{proposition}

\demo{Proof}
Suppose $\beta(M,\Omega_N)>0$ and let $0<b<\beta(M,\Omega_N)$, so \eqref{om1a} is valid for $b$, and then there exists some $q>1$ and $C\ge 1$ such that for all $t\ge 0$ we get $\varphi^{*}_{\omega_{N}}(tq)\ge q\varphi^{*}_{\omega_{M}}(t)+tq\log(q^b)-C$. We apply the Young conjugate to this inequality. Hence for all $s\ge\log(q^b)$ we have
\begin{align*}
\omega_{N}(e^s)&=\varphi_{\omega_{N}}(s)=\varphi^{**}_{\omega_{N}}(s)=\sup_{t\ge 0}\{st-\varphi^{*}_{\omega_{N}}(t)\}
\\&
=\sup_{t\ge 0}\{sqt-\varphi^{*}_{\omega_{N}}(tq)\}\le\sup_{t\ge 0}\{sqt-q\varphi^{*}_{\omega_{M}}(t)-tq\log(q^b)\}+C
\\&
=q\sup_{t\ge 0}\{t(s-\log(q^b))-\varphi^{*}_{\omega_{M}}(t)\}+C=q\varphi^{**}_{\omega_{M}}(s-\log(q^b))+C
\\&
=q\varphi_{\omega_{M}}(s-\log(q^b))+C=q\omega_{M}(e^s/q^b)+C.
\end{align*}
Consequently there does exist $C_1\ge 1$ such that for all $t\ge 0$ we have $\omega_{N}(tq^b)\le q\omega_{M}(t)+C_1$, i.e.  $\limsup_{t\rightarrow+\infty}\frac{\omega_{N}(Kt)}{\omega_{M}(t)}\le K^{1/b}$ with $K:=q^b>1$ and so $(\omega_{M^x},\omega_{M^y})_{\omega_1,1/b'}$ holds with this choice $K$ for any $b'<b$.
Hence we have shown $\alpha(\omega_{M},\omega_{N})\le 1/b'$, and by making $b'$ tend to $b$ and then $b$ tend to $\beta(M,\Omega_N)$, we deduce that  $\alpha(\omega_{M},\omega_{N})\le 1/\beta(M,\Omega_N)$.\vspace{6pt}

Conversely, suppose that $\alpha(\omega_{M},\omega_{N})<\infty$ and let $a>(\omega_{M},\omega_{N})_{\omega_1}$, then $(\omega_{M},\omega_{N})_{\omega_1,a}$ is valid with some $K>1$. So there exist $C\ge 1$ and $0<b<a$ such that $\omega_{N}(Kt)\le K^{b}\omega_{M}(t)+C$ for all $t\ge 0$. Hence by setting $k:=\log(K)>0$ and $s:=\log(t)$ we have $\varphi_{\omega_N}(k+s)=\omega_{N}(e^{k+s})=\omega_N(Kt)\le K^{b}\omega_{M}(e^s)+C=K^{b}\varphi_{\omega_{M}}(s)+C$ for all $s\in\RR$. Applying the Young-conjugate yields for all $s\ge 0$:
\begin{align*}
\varphi^{*}_{\omega_{N}}(s)&=\sup_{t\in\RR}\{st-\varphi_{\omega_{N}}(t)\}=\sup_{t\in\RR}\{s(k+t)-\varphi_{\omega_{N}}(k+t)\}
\\&
\ge\sup_{t\in\RR}\{s(k+t)-K^{b}\varphi_{\omega_{M}}(t)\}-C=\sup_{t\ge 0}\{st-K^{b}\varphi_{\omega_{M}}(t)\}-C+sk
\\&
=K^{b}\sup_{t\ge 0}\{(s/K^{b})t-\varphi_{\omega_{M}}(t)\}-C+sk= K^{b}\varphi^{*}_{\omega_{M}}(s/K^{b})-C+sk.
\end{align*}
Hence we have shown
$$\exists\;K>1\;\exists\;C\ge 1\;\forall\;s'\ge 0:\;\;\;K^{b}\varphi^{*}_{\omega_{M}}(s')+s'K^{b}\log(K)\le \varphi^{*}_{\omega_{N}}(sK^{b})+C,$$
i.e. $$\liminf_{s'\rightarrow+\infty} \frac{\exp\left(\frac{1}{s'K^{b}}\varphi^{*}_{\omega_{N}}(s'K^{b})\right)}{\exp\left(\frac{1}{s'}\varphi^{*}_{\omega_{M}}(s')\right)}\ge\exp(\log(K))=K>(K^b)^{1/a}.
$$
So we have verified $(M,\Omega_N)_{L,1/a}$ with the choice $q:=K^{b}>1$, hence $\beta(M,\Omega_N)\ge 1/a$. When $a$ tends to $(\omega_{M},\omega_{N})_{\omega_1}$, it follows that $\beta(M,\Omega_N)\ge 1/(\omega_{M},\omega_{N})_{\omega_1}$.

One may easily conclude that the stated equality holds in any case (with the conventions $1/0=\infty$, $1/\infty=0$), even if one of the indices is zero or infinity.
\qed\enddemo

Thus, by involving the notation of mixed indices in this section and Proposition \ref{Lomega1comparison}, we can now reformulate Theorem \ref{omega1theorem} as follows.

\begin{theorem}\label{omega1importequivalence}
Let $\mathcal{M}=\{M^x: x\in\mathcal{I}\}$ be a \hyperlink{Msc}{$(\mathcal{M}_{\on{sc}})$} weight matrix and let $\omega_{\mathcal{M}}:=\{\omega_{M^x}: x>0\}$ be the corresponding matrix of associated weight functions.

Then the following are equivalent:
\begin{itemize}
\item[$(i)$] 
Any of the equivalent mixed \hyperlink{om1}{$(\omega_1)$}-conditions of Roumieu type in $(I)$ of Theorem~\ref{omega1theorem} hold true,
\item[$(ii)$]
$$\forall\;x\in\mathcal{I}\;\exists\;y\in\mathcal{I}:\;\;\; \alpha(\omega_{M^x},\omega_{M^y})<\infty,$$

\item[$(iii)$]
$$\forall\;x\in\mathcal{I}\;\exists\;y\in\mathcal{I}:\;\;\; \beta(M^x,\Omega_{M^y})>0.$$
\end{itemize}

Analogously, the following are equivalent:
\begin{itemize}
\item[$(i)$] Any of the equivalent mixed \hyperlink{om1}{$(\omega_1)$}-conditions of Beurling type in $(II)$ of Theorem~\ref{omega1theorem} hold true,
\item[$(ii)$]
$$\forall\;x\in\mathcal{I}\;\exists\;y\in\mathcal{I}:\;\;\; \alpha(\omega_{M^y},\omega_{M^x})<\infty,$$

\item[$(iii)$]
$$\forall\;x\in\mathcal{I}\;\exists\;y\in\mathcal{I}:\;\;\; \beta(M^y,\Omega_{M^x})>0.$$
\end{itemize}
\end{theorem}

\demo{Proof}
Again we limit ourselves to the Roumieu case.

$(i)\Rightarrow(ii)$ It suffices to take into account~\eqref{R-L-consequequstrong}, which easily shows that for every $x$ there exists $y$ such that $(\omega_{M^x},\omega_{M^y})_{\omega_1,a}$ holds for a suitable value of $a$.

$(ii)\Leftrightarrow(iii)$ This is clear from Proposition~\ref{Lomega1comparison}.

$(ii)\Rightarrow(i)$ By hypothesis, for every $x\in\mathcal{I}$ there exists $y\in\mathcal{I}$ and $K>1$ such that
$$
\limsup_{t\rightarrow+\infty}\frac{\omega_{M^y}(Kt)}{\omega_{M^x}(t)}<+\infty.
$$
After a finite iteration (if necessary), we can guarantee that \eqref{R-L-consequ} holds, and we are done.
\qed\enddemo

We close this section with the following observations for given $\mathcal{M}=\{M^x: x\in\mathcal{I}\}$ being \hyperlink{Msc}{$(\mathcal{M}_{\on{sc}})$}:

\begin{itemize}
\item[$(a)$] Lemmas \ref{Lrewriting} and \ref{indexinfinity} imply that condition \hyperlink{R-L}{$(\mathcal{M}_{\{\on{L}\}})$}, respectively  \hyperlink{B-L}{$(\mathcal{M}_{(\on{L})})$}, yields that for every $x\in\mathcal{I}$ there exists $y\in\mathcal{I}$ such that
$\alpha(\omega_{M^x},\omega_{M^y})=0$ and
$\beta(M^x,\Omega_{M^y})=\infty$, resp. $\alpha(\omega_{M^y},\omega_{M^x})=0$ and
$\beta(M^y,\Omega_{M^x})=\infty$.

\item[$(b)$] In particular, if $\Omega$ is the matrix associated with some $\omega\in\hyperlink{omset0}{\mathcal{W}_0}$, by Corollary \ref{Lrewritingcor} we know that, as soon as $\omega$ has in addition \hyperlink{om1}{$(\omega_1)$}, we are in the situation described in $(a)$.
\end{itemize}


\section{Mixed moderate growth conditions}\label{mixedmgsection}
In this section we study the mixed moderate growth conditions for a given weight matrix. The first three equivalent conditions in the next result are stated  in \cite{testfunctioncharacterization} and \cite{dissertation}, new conditions $(iv)$ and $(v)$ deal with replacing the constant $2$ in $(iii)$ by some $q>0$.

\begin{proposition}\label{firstmg}
Let $\mathcal{M}=\{M^x: x\in\mathcal{I}\}$ be \hyperlink{Msc}{$(\mathcal{M}_{\on{sc}})$}, then in the Roumieu setting the following conditions are equivalent:
\begin{itemize}

\item[$(i)$] $\mathcal{M}$ has $\hyperlink{R-mg}{(\mathcal{M}_{\{\on{mg}\}})}$,

\item[$(ii)$] $\forall\;x\in\mathcal{I}\;\exists\;H\ge 1\;\exists\;y\in\mathcal{I}\;\forall\;t\ge 0:\;\;\; 2\omega_{M^y}(t)\le\omega_{M^x}(Ht)+H$,

\item[$(iii)$] $\forall\;x\in\mathcal{I}\;\exists\;C>0\;\exists\;y\in\mathcal{I}\;\forall\;j\in\NN:\;\;\; M^x_{2j}\le C^{2j}(M^y_j)^2$.

 \item[$(iv)$] we have that
 \begin{equation}\label{mgnewroum}
 \exists\;q>0\;\forall\;x\in\mathcal{I}\;\exists\;y\in\mathcal{I}:\;\;\;\limsup_{j\in\NN_{>0},j\rightarrow+\infty}\frac{\exp(\frac{1}{jq}\varphi^{*}_{\omega_{M^x}}(jq))}{\exp(\frac{1}{j}\varphi^{*}_{\omega_{M^y}}(j))}<+\infty,
 \end{equation}

\item[$(v)$] we have that
 \begin{equation}\label{mgnewroumreal}
 \exists\;q>0\;\forall\;x\in\mathcal{I}\;\exists\;y\in\mathcal{I}:\;\;\;\limsup_{j\in\RR_{>0},j\rightarrow+\infty}\frac{\exp(\frac{1}{jq}\varphi^{*}_{\omega_{M^x}}(jq))}{\exp(\frac{1}{j}\varphi^{*}_{\omega_{M^y}}(j))}<+\infty.
 \end{equation}

\end{itemize}

Moreover, in the Beurling setting, we have the following equivalences:
\begin{itemize}

\item[$(i)$] $\mathcal{M}$ has $\hyperlink{B-mg}{(\mathcal{M}_{(\on{mg})})}$,

\item[$(ii)$] $\forall\;x\in\mathcal{I}\;\exists\;H\ge 1\;\exists\;y\in\mathcal{I}\;\forall\;t\ge 0:\;\;\; 2\omega_{M^x}(t)\le\omega_{M^y}(Ht)+H$,

\item[$(iii)$] $\forall\;x\in\mathcal{I}\;\exists\;C>0\;\exists\;y\in\mathcal{I}\;\forall\;j\in\NN:\;\;\; M^y_{2j}\le C^{2j}(M^x_j)^2$.

\item[$(iv)$] we have that
\begin{equation*}
 \exists\;q>0\;\forall\;x\in\mathcal{I}\;\exists\;y\in\mathcal{I}:\;\;\;\limsup_{j\in\NN_{>0},j\rightarrow+\infty}\frac{\exp(\frac{1}{jq}\varphi^{*}_{\omega_{M^y}}(jq))}{\exp(\frac{1}{j}\varphi^{*}_{\omega_{M^x}}(j))}<+\infty,
\end{equation*}

 \item[$(v)$] we have that
 \begin{equation*}
 \exists\;q>0\;\forall\;x\in\mathcal{I}\;\exists\;y\in\mathcal{I}:\;\;\;\limsup_{j\in\RR_{>0},j\rightarrow+\infty}\frac{\exp(\frac{1}{jq}\varphi^{*}_{\omega_{M^y}}(jq))}{\exp(\frac{1}{j}\varphi^{*}_{\omega_{M^x}}(j))}<+\infty.
 \end{equation*}

\end{itemize}
Note: Even if all associated weight functions are equivalent w.r.t. $\hyperlink{sim}{\sim}$, then in general we cannot conclude that \hyperlink{om6}{$(\omega_6)$} for each/some $\omega_{M^x}$ is following.
\end{proposition}

\demo{Proof}
We only treat the Roumieu setting in detail, the Beurling case follows analogously.\vspace{6pt}

First, by combining \cite[Prop. 3.6]{testfunctioncharacterization} and \cite[Thm. 9.5.2, Thm. 9.5.3]{dissertation}, we see that $(i)\Leftrightarrow (ii) \Leftrightarrow (iii)$.\vspace{6pt}

$(iii)\Rightarrow(iv)$ Recall that since each $M^x\in\hyperlink{LCset}{\mathcal{LC}}$, by \eqref{Mxonestable}, we have
$$\forall\;j\in\NN:\;\;\;M^x_j=M^{x;1}_j=\exp(\varphi^{*}_{\omega_{M^x}}(j)).$$
Then $(iii)$ implies the fact that for all $x\in\mathcal{I}$ we find $y\in\mathcal{I}$ such that the \eqref{mgnewroum} is valid with the same indices $x$ and $y$ and the universal choice $q=2$.\vspace{6pt}

$(iv)\Rightarrow(v)$ For any $j\in\RR$, $j\ge 1$, we get $$\frac{\exp(\frac{1}{jq}\varphi^{*}_{\omega_{M^x}}(jq))}{\exp(\frac{1}{j}\varphi^{*}_{\omega_{M^y}}(j))}\le\frac{\exp(\frac{1}{\lfloor j\rfloor 2q}\varphi^{*}_{\omega_{M^x}}(\lfloor j\rfloor 2q))}{\exp(\frac{1}{\lfloor j\rfloor}\varphi^{*}_{\omega_{M^y}}(\lfloor j\rfloor))},$$
since $q\mapsto\frac{1}{jq}\varphi^{*}_{\omega_{M^x}}(jq))$ is nondecreasing for any $x\in\mathcal{I}$ and $j>0$ (fixed) and since $j\le \lfloor j\rfloor+1\le 2\lfloor j\rfloor$ for all $j\ge 1$. Thus \eqref{mgnewroum} with the choice $q>0$ does imply \eqref{mgnewroumreal} with the same indices $x$ and $y$ and with the parameter $q':=q/2$.\vspace{6pt}

$(v)\Rightarrow(iii)$ We have to distinguish between two cases. First, if $q\ge 2$ then we immediately get \eqref{mgnewroum} with $q=2$ ($q\mapsto\frac{1}{jq}\varphi^{*}_{\omega_{M^x}}(jq)$ is nondecreasing for any $x\in\mathcal{I}$ and $j>0$ fixed). Consequently $(i)$ follows again by recalling  \eqref{Mxonestable}.

Second,  if $0<q<2$, then we iterate \eqref{mgnewroumreal} as follows for all $j\ge 1$
$$\frac{\exp(\frac{1}{j2q}\varphi^{*}_{\omega_{M^x}}(j2q))}{\exp(\frac{1}{j}\varphi^{*}_{\omega_{M^z}}(j))}=\frac{\exp(\frac{1}{2jq}\varphi^{*}_{\omega_{M^x}}(2jq))}{\exp(\frac{1}{qj}\varphi^{*}_{\omega_{M^y}}(qj))}\frac{\exp(\frac{1}{jq}\varphi^{*}_{\omega_{M^y}}(jq))}{\exp(\frac{1}{j}\varphi^{*}_{\omega_{M^z}}(j))}.$$
Hence, after applying sufficiently many iterations again depending only on given $q$ (choose $n\in\NN$ minimal to have $2^n q\ge 2$), we get \eqref{mgnewroumreal} with $q'\ge 2$ and some index $z\in\mathcal{I}$ and are again able to conclude.
\qed\enddemo

Assertion $(ii)$ in the previous result is the mixed \hyperlink{om6}{$(\omega_6)$}-condition of the particular, Roumieu or Beurling, type. Using resp. applying iterations as in the proofs of the previous section it is straight-forward to extend the list(s) of equivalences in Proposition \ref{firstmg} by replacing in $(ii)$ the value $2$ by any $C>1$ not depending on $x$ and $y$.

\subsection{Mixed growth indices based on mixed moderate growth conditions}\label{mixedmgindexsection}
Let $M,N\in\hyperlink{LCset}{\mathcal{LC}}$ with $M\le N$ be given and $a>0$. We write $(N,\Omega_M)_{\on{mg},a}$, if
\begin{equation}\label{mga}
\exists\;q>1:\;\limsup_{j\rightarrow+\infty}\frac{\exp(\frac{1}{jq}\varphi^{*}_{\omega_{M}}(jq))}{\exp(\frac{1}{j}\varphi^{*}_{\omega_{N}}(j))}<q^{a},
\end{equation}
which should be compared with \cite[Thm. 3.16 $(v)$]{index}. Note that, if we write $\hyperlink{OmomM}{\Omega_M}=\{W^q=(W_j^q)_{j\in\NN}\colon q>0\}$ for the weight matrix associated with $\omega_M$, so in particular $W^1=M$ (see Section \ref{weightmatrixfromfunction}), then the previous condition reads
$$
\exists\;q>1:\;\;\;\limsup_{j\rightarrow+\infty} \left(\frac{W_j^q}{N_j}\right)^{1/j} <q^{a}.
$$
Given weight two functions $\omega$ and $\sigma$ with $\sigma\ge\omega$ we write $(\omega,\sigma)_{\omega_6,b}$ if
\begin{equation}\label{mg1a}
\exists\;K>1:\;\liminf_{t\rightarrow+\infty}\frac{\sigma(Kt)}{\omega(t)}>K^b,
\end{equation}
which should be compared with \cite[Thm. 2.16 $(iii)$]{index}.
It is immediate to see that if \eqref{mga} is satisfied for some $a>0$, then also for all $a'>a$ with the same choices of $q$ and if \eqref{mg1a} is satisfied for some given $b>0$, then also for all $0<b'<b$ with the same choices for $K$. According to these growth restrictions we can put
\begin{equation}\label{mg2a}
\alpha(N,\Omega_M):=\inf\{a>0: (N,\Omega_M)_{\on{mg},a}\}\in [0,\infty],
\end{equation}
and
\begin{equation}\label{mg3a}
\beta(\omega,\sigma):=\sup\{b>0: (\omega,\sigma)_{\omega_6,b}\}\in [0,\infty].
\end{equation}

If the corresponding sets are empty, we write $\alpha(N,\Omega_M)=\infty$ and $\beta(\omega,\sigma)=0$.

The next result shows that the values defined in \eqref{mg2a} and \eqref{mg3a} are related.

\begin{proposition}\label{mgomega6comparison}
Let $M,N\in\hyperlink{LCset}{\mathcal{LC}}$ be given with $M\le N$ and $\omega_{M}, \omega_{N}$ the corresponding associated weight functions. Then we get
$$\alpha(N,\Omega_M)=\frac{1}{\beta(\omega_{N},\omega_{M})}.$$
\end{proposition}

\demo{Proof}
Let  $a>\alpha(N,\Omega_M)$, so \eqref{mga} is valid for $a>0$. Then there exists some $C\ge 1$ such that for all $t\ge 0$ we get $\varphi^{*}_{\omega_{M}}(tq)\le q\varphi^{*}_{\omega_{N}}(t)+tq\log(q^{1/a})+C$. We apply the Young conjugate to this inequality, hence for all $s\ge\log(q^{a})$ we have
\begin{align*}
\omega_{M}(e^s)&=\varphi_{\omega_{M}}(s)=\varphi^{**}_{\omega_{M}}(s)=\sup_{t\ge 0}\{st-\varphi^{*}_{\omega_{M}}(t)\}
\\&
=\sup_{t\ge 0}\{sqt-\varphi^{*}_{\omega_{M}}(tq)\}\ge\sup_{t\ge 0}\{sqt-q\varphi^{*}_{\omega_{N}}(t)-tq\log(q^{a})\}-C
\\&
=q\sup_{t\ge 0}\{t(s-\log(q^{a}))-\varphi^{*}_{\omega_{N}}(t)\}-C=q\varphi^{**}_{\omega_{N}}(s-\log(q^{a}))-C
\\&
=q\varphi_{\omega_{N}}(s-\log(q^{a}))-C=q\omega_{N}(e^s/q^{a})-C.
\end{align*}
Consequently, there does exist $C_1\ge 1$ such that for all $t\ge 0$ we have $q\omega_{N}(t)\le\omega_{M}(tq^{a})+C_1$, i.e. $\liminf_{t\rightarrow+\infty}\frac{\omega_{M}(Kt)}{\omega_{N}(t)}\ge K^{1/a}$ with $K:=q^{a}>1$. So $(\omega_{N},\omega_{M})_{\omega_6,1/a'}$ is valid with this choice $K$ for any $a'>a$ which proves $\beta(\omega_{N},\omega_{M})\ge 1/a'$. Since $a$ can be chosen arbitrarily close to $\alpha(N,\Omega_M)$  we have shown $\beta(\omega_{N},\omega_{M})\ge 1/\alpha(N,\Omega_M)$.\vspace{6pt}

Conversely, let now $0<b<\beta(\omega_{N},\omega_{M})$, i.e. $(\omega_{N},\omega_{M})_{\omega_6,b}$ is valid with some $K>1$, so there exists $C\ge 1$  such that $\omega_{M}(Kt)\ge K^b\omega_{N}(t)-C$ for all $t\ge 0$. By setting $k:=\log(K)>0$ and $s:=\log(t)$ we have $\varphi_{\omega_{M}}(k+s)=\omega_{M}(e^{k+s})\ge K^b\omega_{N}(e^s)-C=K^b\varphi_{\omega_{N}}(s)-C$ for all $s\in\RR$. Hence applying the Young-conjugate yields for all $s\ge 0$:
\begin{align*}
\varphi^{*}_{\omega_{M}}(s)&=\sup_{t\in\RR}\{st-\varphi_{\omega_{M}}(t)\}=\sup_{t\in\RR}\{s(k+t)-\varphi_{\omega_{M}}(k+t)\}
\\&
\le\sup_{t\in\RR}\{s(k+t)-K^b\varphi_{\omega_{N}}(t)\}+C=\sup_{t\ge 0}\{st-K^b\varphi_{\omega_{N}}(t)\}+C+sk
\\&
=K^b\sup_{t\ge 0}\{(s/K^b)t-\varphi_{\omega_{N}}(t)\}+C+sk=K^b\varphi^{*}_{\omega_{N}}(s/K^b)+C+sk.
\end{align*}
Thus we have shown that there exists $C\geq 1$ such that $\varphi^{*}_{\omega_{M}}(tK^b)\le K^b\varphi^{*}_{\omega_{N}}(t)+tK^b\log(K)+C $ for all $t\ge 0$.
Consequently,  $(N,\Omega_{M})_{\omega_6,1/b'}$ is valid with the choice $q:=K^b$ for any $0<b'<b$ which proves $1/b'\ge \alpha(N,\Omega_{M})$. Hence  we obtain $\beta(\omega_{N},\omega_{M}) \le 1/\alpha(N, \Omega_M)$.\vspace{6pt}

Note that this proof also shows that if one of the indices is zero the other is infinity, so also in this situation the equality holds with the usual convention.
\qed\enddemo

\begin{remark}
By inspecting the proofs of Propositions \ref{Lomega1comparison} and \ref{mgomega6comparison} we see that the convexity assumption \hyperlink{om4}{$(\omega_4)$} is indispensable in order to ensure $\varphi^{**}_{\omega}=\varphi_{\omega}$. Note that this condition is always satisfied for any $\omega_M$ when having $M\in\hyperlink{LCset}{\mathcal{LC}}$, see Lemma \ref{assoweightomega0}.
\end{remark}


Thus by Proposition \ref{firstmg} and Proposition \ref{mgomega6comparison}, we immediately get the following characterization.


\begin{theorem}\label{mixedmgcharact}
Let $\mathcal{M}=\{M^x: x>0\}$ be a \hyperlink{Msc}{$(\mathcal{M}_{\on{sc}})$} weight matrix and $\omega_{\mathcal{M}}:=\{\omega_{M^x}: x>0\}$ be the corresponding matrix of associated weight functions. Then the following conditions are equivalent:
\begin{itemize}
\item[$(i)$] $\hyperlink{R-mg}{(\mathcal{M}_{\{\on{mg}\}})}$ holds true,

\item[$(ii)$] $\forall\;x\in\mathcal{I}\;\exists\;y\in\mathcal{I}:\;\;\;\alpha(M^y,\Omega_{M^x}) <\infty,$

\item[$(iii)$] $\forall\;x\in\mathcal{I}\;\exists\;y\in\mathcal{I}:\;\;\;\beta(\omega_{M^y},\omega_{M^x})>0.$
\end{itemize}

Analogously, for the Beurling case, the following are equivalent:
\begin{itemize}
\item[$(i)$] $\hyperlink{B-mg}{(\mathcal{M}_{(\on{mg})})}$ holds true,

\item[$(ii)$] $\forall\;x\in\mathcal{I}\;\exists\;y\in\mathcal{I}:\;\;\;\alpha(M^x,\Omega_{M^y}) <\infty,$

\item[$(iii)$] $\forall\;x\in\mathcal{I}\;\exists\;y\in\mathcal{I}:\;\;\;\beta(\omega_{M^x},\omega_{M^y})>0.$
\end{itemize}
\end{theorem}

\section{Consequences for ultradifferentiable classes}\label{consequencesforultra}
The aim of this final section is to extend the characterizing results from \cite{BonetMeiseMelikhov07} to the matrix setting, to refine the consequences shown in \cite{testfunctioncharacterization} and to give the new mixed growth indices introduced in the previous sections an interpretation in terms of the characterization for the equivalence of ultradifferentiable classes defined by (abstractly given) weight matrices.

\subsection{Classes defined by abstractly given weight matrices}\label{consequencesforultrasequ}
First let us recall the main result \cite[Theorem 3.2]{testfunctioncharacterization}.

\begin{theorem}\label{Thm32testfunc}
Let $\mathcal{M}$ be a \hyperlink{Msc}{$(\mathcal{M}_{\on{sc}})$} weight matrix and $\omega_{\mathcal{M}}$ be the corresponding matrix of associated weight functions. If $\mathcal{M}$ does have $(\mathcal{M}_{[\on{mg}]})$ and $(\mathcal{M}_{[\on{L}]})$, then
\begin{equation}\label{Thm32testfuncequ}
\mathcal{E}_{[\mathcal{M}]}=\mathcal{E}_{[\omega_{\mathcal{M}}]}
\end{equation}
and the equality holds as top. vector spaces.
\end{theorem}

The aim of this section is to show that \eqref{Thm32testfuncequ} is characterized in terms of the growth properties studied before. More precisely we will see that assumption $(\mathcal{M}_{[\on{L}]})$ is too strong.

The proof of Theorem \ref{Thm32testfunc} has been split into two parts. The first one \cite[Theorem 3.4]{testfunctioncharacterization} deals with the mixed moderate growth conditions and can be reformulated as follows.

\begin{theorem}\label{Thm34testfunc}
Let $\mathcal{M}:=\{M^x: x\in\mathcal{I}\}$ be a \hyperlink{Msc}{$(\mathcal{M}_{\on{sc}})$} weight matrix and recall $\mathcal{M}^{(2)}:=\{M^{x;l}: x\in\mathcal{I}, l>0\}$, then the following conditions are equivalent:
\begin{itemize}
\item[$(i)$] One has
$$\mathcal{E}_{\{\mathcal{M}\}}=\mathcal{E}_{\{\mathcal{M}^{(2)}\}},$$
as top. vector spaces,

\item[$(ii)$] we have that
$$\forall\;x\in\mathcal{I}\;\exists\;y\in\mathcal{I}:\;\;\;\alpha(M^y,\Omega_{M^x}) <\infty,$$

\item[$(iii)$] we have that
$$\forall\;x\in\mathcal{I}\;\exists\;y\in\mathcal{I}:\;\;\;\beta(\omega_{M^y},\omega_{M^x})>0.$$

\end{itemize}

Analogously,  the following conditions are equivalent:
\begin{itemize}
\item[$(i)$] One has
$$\mathcal{E}_{(\mathcal{M})}=\mathcal{E}_{(\mathcal{M}^{(2)})},$$
as top. vector spaces,

\item[$(ii)$] we have that
$$\forall\;x\in\mathcal{I}\;\exists\;y\in\mathcal{I}:\;\;\;\alpha(M^x,\Omega_{M^y}) <\infty,$$

\item[$(iii)$] we have that
$$\forall\;x\in\mathcal{I}\;\exists\;y\in\mathcal{I}:\;\;\;\beta(\omega_{M^x},\omega_{M^y})>0.$$
\end{itemize}
\end{theorem}

\demo{Proof} We prove only the statement in the Roumieu case.

The equivalence of $(i)$ and condition $\hyperlink{R-mg}{(\mathcal{M}_{\{\on{mg}\}})}$ is shown in \cite[Theorem 3.7, Proposition 3.9]{testfunctioncharacterization}.

Theorem \ref{mixedmgcharact} implies the equivalence of the three statements.
\qed\enddemo

The second part has been treated in \cite[Section 3.10]{testfunctioncharacterization} by studying consequences of the assumption $(\mathcal{M}_{[\on{L}]})$ on $\mathcal{M}$. However, in order to have equality between the classes $\mathcal{E}_{[\mathcal{M}^{(2)}]}$ and $\mathcal{E}_{[\omega_{\mathcal{M}}]}$ it is sufficient to have one of the equivalent but weaker conditions from Theorem \ref{omega1theorem} above.

\begin{theorem}\label{Thm311testfunc}
Let $\mathcal{M}:=\{M^x: x\in\mathcal{I}\}$ be a \hyperlink{Msc}{$(\mathcal{M}_{\on{sc}})$} weight matrix, let $\omega_{\mathcal{M}}$ be the corresponding matrix of associated weight functions and finally $\mathcal{M}^{(2)}:=\{M^{x;l}: x\in\mathcal{I}, l>0\}$.

The following assertions are equivalent:
\begin{itemize}
\item[$(i)$] One has
$$\mathcal{E}_{\{\mathcal{M}^{(2)}\}}=\mathcal{E}_{\{\omega_{\mathcal{M}}\}}$$
as top. vector spaces,

\item[$(ii)$] we have that
$$\forall\;x\in\mathcal{I}\;\exists\;y\in\mathcal{I}:\;\;\; \alpha(\omega_{M^x},\omega_{M^y})<\infty,$$

\item[$(iii)$] we have that
$$\forall\;x\in\mathcal{I}\;\exists\;y\in\mathcal{I}:\;\;\; \beta(M^x,\Omega_{M^y})>0.$$
\end{itemize}

Analogously, the following are equivalent:
\begin{itemize}
\item[$(i)$] One has $$\mathcal{E}_{(\mathcal{M}^{(2)})}=\mathcal{E}_{(\omega_{\mathcal{M}})}$$
as top. vector spaces,

\item[$(ii)$] we have that
$$\forall\;x\in\mathcal{I}\;\exists\;y\in\mathcal{I}:\;\;\; \alpha(\omega_{M^y},\omega_{M^x})<\infty,$$

\item[$(iii)$] we have that
$$\forall\;x\in\mathcal{I}\;\exists\;y\in\mathcal{I}:\;\;\; \beta(M^y,\Omega_{M^x})>0.$$
\end{itemize}
\end{theorem}


\demo{Proof}
In the Roumieu case, $(ii)$ and $(iii)$ are equivalent, by Theorem \ref{omega1importequivalence}, to the fact that one/each of the conditions from $(I)$ in Theorem \ref{omega1theorem} holds true. These conditions are shown to imply $(i)$ in \cite[Section 3.10]{testfunctioncharacterization}. The same arguments apply for these implications in the Beurling case. So, it is only pending the proof that $(i)$ in the particular Roumieu or Beurling case implies one of the equivalent conditions from $(I)$, resp. $(II)$, in Theorem \ref{omega1theorem}.\vspace{6pt}

{\itshape The Roumieu case.} By assumption the nontrivial inclusion $\mathcal{E}_{\{\mathcal{M}^{(2)}\}}\subseteq\mathcal{E}_{\{\omega_{\mathcal{M}}\}}$ is valid.

Then recall that for any given $N\in\hyperlink{LCset}{\mathcal{LC}}$ we can define (recall $\nu_k:=\frac{N_k}{N_{k-1}}$, $\nu_0:=1$):
$$\theta_{N}(x):=\sum_{k=0}^{+\infty}\frac{N_k}{(2\nu_k)^k}\exp(2i\nu_k
x), \quad x \in \RR.$$
We get that $\theta_N\in\mathcal{E}_{\{N\}}(\RR,\CC)$ (in fact $\theta_N$ does admit global estimates on whole $\RR$) and
\begin{equation*}
\theta_N^{(j)}(0)=i^j s_{j}  \text{ with } s_j\ge N_{j}\, ,  \quad \forall\,j\in\NN.
\end{equation*}
We refer to \cite[Theorem 1]{thilliez}, for a detailed proof see also \cite[Prop. 3.1.2]{diploma} and \cite[Lemma 2.9]{compositionpaper}. However, it is not difficult to see that $\theta_N$ does not belong to the Beurling type class $\mathcal{E}_{(N)}(\RR,\CC)$.

Each $M^{x;a}\in\mathcal{M}^{(2)}$ does belong to the class \hyperlink{LCset}{$\mathcal{LC}$}. Let $h>1$ be arbitrary, but from now on fixed. Let also $x\in\mathcal{I}$ be arbitrary but fixed and put $\widetilde{M}^{x;1}:=(h^jM^{x;1}_j)_{j\in\NN}$. Clearly $\widetilde{M}^{x;1}\in\hyperlink{LCset}{\mathcal{LC}}$ and this sequence is equivalent to $M^{x;1}$. Consequently $\theta_{\widetilde{M}^{x;1}}\in\mathcal{E}_{\{\widetilde{M}^{x;1}\}}(\RR,\CC)=\mathcal{E}_{\{M^{x;1}\}}(\RR,\CC)\subseteq\mathcal{E}_{\{\mathcal{M}^{(2)}\}}(\RR,\CC)$ and by the inclusion $\mathcal{E}_{\{\mathcal{M}^{(2)}\}}(\RR,\CC)\subseteq\mathcal{E}_{\{\omega_{\mathcal{M}}\}}(\RR,\CC)$ we get
$$\exists\;y\in\mathcal{I}\;\exists\;b>0\;\exists\;D>0\;\forall\;j\in\NN:\;\;\;h^jM^{x;1}_j=\widetilde{M}^{x;1}_j\le|\theta_{\widetilde{M}^{x;1}}(0)|\le D\exp\left(\frac{1}{b}\varphi^{*}_{\omega_{M^y}}(bj)\right)=DM^{y;b}_j,$$
hence one may easily deduce that \eqref{naturalitycor1new1} is verified.\vspace{6pt}

{\itshape The Beurling case.} We follow the ideas from the characterization of the inclusion relations, see \cite[Prop. 4.6]{compositionpaper} and also \cite[Prop. 3.9]{testfunctioncharacterization}. By assumption the nontrivial inclusion $\mathcal{E}_{(\omega_{\mathcal{M}})}(\RR,\CC)\subseteq\mathcal{E}_{(\mathcal{M}^{(2)})}(\RR,\CC)$ is valid and the inclusion mapping is continuous by the closed graph theorem. Note that both spaces are Fr\'{e}chet spaces.\vspace{6pt}

This means that
$\forall\;K\subseteq\RR\;\text{compact}\;\forall\;x\in\mathcal{I}\;\forall\;a>0\;\forall\;h>0\;\exists\;K_1\subseteq\RR\;\text{compact}\;\exists\;y\in\mathcal{I}\;\exists\;b>0\;\exists\;D>0\;\forall\;f\in\mathcal{E}_{(\omega_{\mathcal{M}})}(\RR,\CC):\;\;\;\|f\|_{M^{x;a},K,h}\le D\|f\|_{\omega_{M^y},K_1,b}$.

We apply this to the functions $f_s(t):=e^{its}$, $s\ge 0$ and $t\in\RR$. Note that $f_s\in\mathcal{E}_{(\omega_{\mathcal{M}})}(\RR,\CC)$ because $|f^{(j)}_s(t)|=s^j$ for all $j\in\NN$, $s\ge 0$ and $t\in\RR$ and $(M^{x;a}_j)^{1/j}\rightarrow+\infty$ as $j\rightarrow+\infty$ for each $x\in\mathcal{I}$ and $a>0$. Then we estimate as follows:
\begin{align*}
\sup_{j\in\NN}\frac{s^j}{h^jM^{x;a}_j}&=\sup_{j\in\NN, t\in K}\frac{|f_s^{(j)}(t)|}{h^jM^{x;a}_j}=\|f_s\|_{M^{x;a},K,h}\le D\|f_s\|_{\omega_{M^y},K_1,b}=D\sup_{j\in\NN, t\in K_1}\frac{|f_s^{(j)}(t)|}{\exp(\frac{1}{b}\varphi^{*}_{\omega_{M^y}}(bj))}
\\&
=D\sup_{j\in\NN}\frac{s^j}{M^{y;b}_j},
\end{align*}
which yields $\exp(\omega_{M^{x;a}}(s/h))\le D\exp(\omega_{M^{y;b}}(s))$ for all $s\ge 0$. Hence, by applying \eqref{Prop32Komatsu} we get for all $j\in\NN$
$$M^{y;b}_j=\sup_{t\ge 0}\frac{t^j}{\exp(\omega_{M^{y;b}}(t))}\le D\sup_{t\ge 0}\frac{t^j}{\exp(\omega_{M^{x;a}}(t/h))}=Dh^jM^{x;a}_j.$$
Consequently we have shown
$$\forall\;x\in\mathcal{I}\;\forall\;a>0\;\forall\;C> 1\;\exists\;y\in\mathcal{I}\;\exists\;b>0\;\exists\;D>0\;\forall\;j\in\NN:\;\;\;C^jM^{y;b}_j\le DM^{x;a}_j,$$
i.e. \eqref{naturalitycor1roumnewnewbeur} and are done.
\qed\enddemo

Gathering all this information we are able to formulate the following characterization.

\begin{theorem}\label{Thm32testfuncnew}
Let $\mathcal{M}$ be a \hyperlink{Msc}{$(\mathcal{M}_{\on{sc}})$} weight matrix and $\omega_{\mathcal{M}}$ be the corresponding matrix of associated weight functions. Then the following are equivalent:

\begin{itemize}
\item[$(i)$] $\mathcal{M}$ does have $\hyperlink{R-mg}{(\mathcal{M}_{\{\on{mg}\}})}$ and one/each of the characterizing conditions from $(I)$ in Theorem \ref{omega1theorem} (e.g. the mixed \hyperlink{om1}{$(\omega_1)$}-conditions of Roumieu type),

\item[$(ii)$] we have as top. vector spaces the equalities
$$\mathcal{E}_{\{\mathcal{M}\}}=\mathcal{E}_{\{\mathcal{M}^{(2)}\}}=\mathcal{E}_{\{\omega_{\mathcal{M}}\}},$$

\item[$(iii)$] $\forall\;x\in\mathcal{I}\;\exists\;y\in\mathcal{I}:\;\;\;\alpha(M^y,\Omega_{M^x}) <\infty,$ and\par
$\forall\;x\in\mathcal{I}\;\exists\;y\in\mathcal{I}:\;\;\; \beta(M^x,\Omega_{M^y})>0.$

\item[$(iv)$] $\forall\;x\in\mathcal{I}\;\exists\;y\in\mathcal{I}:\;\;\;\beta(\omega_{M^y},\omega_{M^x})>0$, and\par
$\forall\;x\in\mathcal{I}\;\exists\;y\in\mathcal{I}:\;\;\; \alpha(\omega_{M^x},\omega_{M^y})<\infty.$

\end{itemize}

Analogously, the following are equivalent:

\begin{itemize}
\item[$(i)$] $\mathcal{M}$ does have $\hyperlink{B-mg}{(\mathcal{M}_{(\on{mg})})}$ and one/each of the characterizing conditions from $(II)$ in Theorem \ref{omega1theorem} (e.g. the mixed \hyperlink{om1}{$(\omega_1)$}-conditions of Beurling type),

\item[$(ii)$] we have as top. vector spaces the equalities
$$\mathcal{E}_{(\mathcal{M})}=\mathcal{E}_{(\mathcal{M}^{(2)})}=\mathcal{E}_{(\omega_{\mathcal{M}})},$$

\item[$(iii)$] $\forall\;x\in\mathcal{I}\;\exists\;y\in\mathcal{I}:\;\;\; \beta(M^y,\Omega_{M^x})>0$, and\par
$\forall\;x\in\mathcal{I}\;\exists\;y\in\mathcal{I}:\;\;\;\alpha(M^x,\Omega_{M^y}) <\infty,$

\item[$(iv)$]
$\forall\;x\in\mathcal{I}\;\exists\;y\in\mathcal{I}:\;\;\; \alpha(\omega_{M^y},\omega_{M^x})<\infty$, and\par
$\forall\;x\in\mathcal{I}\;\exists\;y\in\mathcal{I}:\;\;\;\beta(\omega_{M^x},\omega_{M^y})>0.$

\end{itemize}
\end{theorem}

We give now some comments on ''classical situations'' to which this main result applies.

\begin{itemize}
\item[$(a)$] Theorem \ref{Thm32testfuncnew} applies to $\mathcal{M}\equiv\Omega=\{W^x: x>0\}$, when $\Omega$ is associated with a given weight function $\omega\in\hyperlink{omset0}{\mathcal{W}_0}$ (see Section \ref{weightmatrixfromfunction}). In this situation this theorem is consistent with the known results from \cite{BonetMeiseMelikhov07} and \cite{compositionpaper}: On the one hand $\Omega$ satisfies all requirements from assertion $(i)$ (for both types) provided that $\omega\in\hyperlink{omset1}{\mathcal{W}}$, see \eqref{newmoderategrowth} and Corollary \ref{Lrewritingcor}. On the other hand, we have in this situation
$$\mathcal{E}_{[\Omega]}=\mathcal{E}_{[\omega]}=\mathcal{E}_{[\omega_{\Omega}]},$$
with the first equality holding by \cite[Theorem 5.14 $(2)$]{compositionpaper} and the second one by \cite[Lemma 5.7]{compositionpaper}. Note that $\omega\hyperlink{sim}{\sim}\omega_{W^x}$ for all $x>0$ (i.e. $(vi)$ in Sect. \ref{weightmatrixfromfunction}) , hence the matrix $\omega_{\Omega}$ is constant.\vspace{6pt}

Alternatively, by our developed knowledge on mixed growth indices we can also directly verify assertions $(iv)$: By \eqref{newmoderategrowth} and Theorem \ref{mixedmgcharact} we get both requirements on the index $\beta$; whereas Corollary \ref{Lrewritingcor} and Lemma \ref{indexinfinity} yield that both requirements on the index $\alpha$ are valid, more precisely even with the value $0$.

\item[$(b)$] If $\mathcal{M}=\{M\}$, then $\omega_{\mathcal{M}}=\{\omega_M\}$ and there is no difference between Roumieu- and Beurling-like conditions. In this case Theorem \ref{Thm32testfuncnew} yields the characterizations shown in \cite{BonetMeiseMelikhov07}. More precisely, we remark that in $(i)$ assumption \hyperlink{mg}{$(\on{mg})$} on $M$ implies that some/any assertions in Theorem \ref{omega1theorem} are valid if and only if
\begin{equation}\label{beta3}
\exists\;Q\in\NN_{\ge 2}:\;\;\;\liminf_{j\rightarrow+\infty}\frac{\mu_{Qj}}{\mu_j}>1,
\end{equation}
which follows by \cite[Prop. 3.4]{subaddlike}. However, note also that in \cite{BonetMeiseMelikhov07} the general assumptions on $M$ have been slightly stronger than in our result, i.e. apart from $M\in\hyperlink{LCset}{\mathcal{LC}}$ also $\liminf_{p\rightarrow\infty}(M_p/p!)^{1/p}>0$ (this is precisely $(M0)$ in \cite{BonetMeiseMelikhov07}) and $(M2')$ from \cite{Komatsu73} ({\itshape derivation closedness}) are required.\vspace{6pt}

Finally, let us see $(i)\Leftrightarrow(iv)$ directly and not using any information on the underlying ultradifferentiable classes: In this situation assertion $(iv)$ reads $\alpha(\omega_M)<\infty$ and $\beta(\omega_M)>0$, which is by \cite[Cor. 2.14, Cor. 2.17]{index} equivalent to requiring both \hyperlink{om1}{$(\omega_1)$} and \hyperlink{om6}{$(\omega_6)$} for $\omega_M$. By combining \cite[Prop. 3.6]{Komatsu73} and \cite[Thm. 3.1]{subaddlike} we obtain the equivalence to $(i)$.
\end{itemize}

We close with the following observation concerning the stability of the arising conditions w.r.t. the natural equivalence relations $[\approx]$ characterizing the permanence of the class $\mathcal{E}_{[\mathcal{M}]}$.

\begin{remark}\label{equivremark}
If both $\mathcal{M}$ and $\mathcal{N}$ are \hyperlink{Msc}{$(\mathcal{M}_{\on{sc}})$}, then in \cite[Prop. 4.6]{compositionpaper} it has been shown that $\mathcal{E}_{[\mathcal{M}]}=\mathcal{E}_{[\mathcal{N}]}$ if and only if $\mathcal{M}[\approx]\mathcal{N}$.\vspace{6pt}

And this is consistent with Theorem \ref{Thm32testfuncnew}: For this recall that $(\mathcal{M}_{[\on{mg}]})$ is clearly stable under relation $[\approx]$. Similarly this holds true for the equivalent conditions from Theorem \ref{omega1theorem} of the particular, Roumieu or Beurling, type; it can be directly checked for $(I)(iv)$ resp. $(II)(iv)$: \eqref{omega1mixedcharactequmodstrong} resp. \eqref{omega1mixedcharactequbeurmodstrong} are clearly preserved under $\{\approx\}$ resp. $(\approx)$.

Recall that in general it is not clear that the stronger condition $(\mathcal{M}_{[\on{L}]})$ is preserved under $[\approx]$.
\end{remark}

\subsection{Classes defined by an abstractly given matrix of weight functions}\label{consequencesforultrafct}
The aim of this section is to apply the main result from the previous section to a weight matrix $\mathcal{N}$ which is obtained by an abstractly given matrix of weight functions $\mathcal{M}_{\mathcal{W}}:=\{\omega^x\in\hyperlink{omset0}{\mathcal{W}_0}: x\in\mathcal{I}\}$. Recall that here $\omega^y\le\omega^x$ for any $x,y\in\mathcal{I}$ with $x\le y$ and the class $\mathcal{E}_{[\mathcal{M}_{\mathcal{W}}]}$ is defined analogously as $\mathcal{E}_{[\omega_{\mathcal{M}}]}$ in Section \ref{classes}.

\begin{theorem}\label{generalweightmatrixthm}
Let $\mathcal{M}_{\mathcal{W}}:=\{\omega^x\in\hyperlink{omset0}{\mathcal{W}_0}: x\in\mathcal{I}\}$ be given and let the \hyperlink{Msc}{$(\mathcal{M}_{\on{sc}})$} weight matrix $\mathcal{N}=\{N^x: x\in\mathcal{I}\}$ be defined by
\begin{equation}\label{generalweightmatrixthmequ}
N^x_p:=\exp(\varphi^{*}_{\omega^x}(p)).
\end{equation}

\begin{itemize}
\item[$(I)$] The following are equivalent in the Roumieu case:

\begin{itemize}
\item[$(i)$] $\mathcal{N}$ does have $\hyperlink{R-mg}{(\mathcal{M}_{\{\on{mg}\}})}$ and one/each of the characterizing conditions from $(I)$ in Theorem \ref{omega1theorem} (e.g. the mixed \hyperlink{om1}{$(\omega_1)$}-conditions of Roumieu type),

\item[$(ii)$] we have as top. vector spaces the equality
$$\mathcal{E}_{\{\mathcal{M}_{\mathcal{W}}\}}=\mathcal{E}_{\{\mathcal{N}\}}.$$


\item[$(iii)$] $\forall\;x\in\mathcal{I}\;\exists\;y\in\mathcal{I}:\;\;\;\beta(\omega^{y},\omega^{x})>0$, and\par
$\forall\;x\in\mathcal{I}\;\exists\;y\in\mathcal{I}:\;\;\; \alpha(\omega^{x},\omega^{y})<\infty.$

\end{itemize}

\item[$(II)$] Analogously, the following are equivalent in the Beurling case:

\begin{itemize}
\item[$(i)$] $\mathcal{N}$ does have $\hyperlink{B-mg}{(\mathcal{M}_{(\on{mg})})}$ and one/each of the characterizing conditions from $(II)$ in Theorem \ref{omega1theorem} (e.g. the mixed \hyperlink{om1}{$(\omega_1)$}-conditions of Beurling type),

\item[$(ii)$] we have as top. vector spaces the equality
$$\mathcal{E}_{(\mathcal{M}_{\mathcal{W}})}=\mathcal{E}_{(\mathcal{N})}.$$


\item[$(iii)$]
$\forall\;x\in\mathcal{I}\;\exists\;y\in\mathcal{I}:\;\;\; \alpha(\omega^{y},\omega^{x})<\infty$, and\par
$\forall\;x\in\mathcal{I}\;\exists\;y\in\mathcal{I}:\;\;\;\beta(\omega^{x},\omega^{y})>0.$

\end{itemize}

\end{itemize}
\end{theorem}

\demo{Proof}

First, by \eqref{generalweightmatrixthmequ} the matrix $\mathcal{N}$ is \hyperlink{Msc}{$(\mathcal{M}_{\on{sc}})$} and by $(vi)$ in Sect. \ref{weightmatrixfromfunction} we obtain
\begin{equation}\label{weighfctthmequ}
\forall\;x\in\mathcal{I}:\;\;\;\omega_{N^x}\hyperlink{sim}{\sim}\omega^x.
\end{equation}
This fact implies (e.g. see the proof of \cite[Lemma 5.16 $(1)$]{compositionpaper}) that $\mathcal{E}_{[\omega_{\mathcal{N}}]}=\mathcal{E}_{[\mathcal{M}_{\mathcal{W}}]}$. Now, we are in a position to apply Theorem \ref{Thm32testfuncnew} to the matrix $\mathcal{N}$, so that
$$\mathcal{E}_{[\mathcal{N}]}=\mathcal{E}_{[\mathcal{N}^{(2)}]}=\mathcal{E}_{[\omega_{\mathcal{N}}]},$$
and the conclusion is straightforward both in the Roumieu and the Beurling case.
\qed\enddemo

If $\mathcal{M}_{\mathcal{W}}:=\{\omega\}$ and so $\mathcal{N}=\{N\}$, then there is no difference between the Roumieu and the Beurling case and so for the mixed indices we precisely get the indices $\alpha(\omega)$ and $\beta(\omega)$ studied in \cite{index}. Consequently, assertion $(iii)$ in Theorem \ref{generalweightmatrixthm} amounts to having both \hyperlink{om1}{$(\omega_1)$} and \hyperlink{om6}{$(\omega_6)$} for $\omega$, see \cite[Cor. 2.14, Cor. 2.17]{index}. So this result is consistent with \cite[Cor. 16]{BonetMeiseMelikhov07}. However, note that there \hyperlink{om1}{$(\omega_1)$} is a standard assumption, i.e. $\omega\in\hyperlink{omset1}{\mathcal{W}}$. Finally, concerning $N$ and assertion $(i)$ we recall: As mentioned in the previous section, since $N$ has \hyperlink{mg}{$(\on{mg})$} we have that $N$ satisfies some/any of the assertions in Theorem \ref{omega1theorem} if and only if \eqref{beta3} holds.

\vspace{6pt}
\textbf{Acknowledgements}: The first two authors are supported by
the Spanish Ministry of Science and Innovation under the project PID2019-105621GB-I00.
The third author is supported by FWF-Project P~33417-N.\par

\bibliographystyle{plain}
\bibliography{Bibliography}

\vskip1cm

\textbf{Affiliation}:\\
J.~Jim\'{e}nez-Garrido:\\
Departamento de Matem\'aticas, Estad{\'\i}stica y Computaci\'on,\\
Universidad de Cantabria,\\
Facultad de Ciencias, Avda. de los Castros, 48, 39005 Santander, Spain.\\

Instituto de Investigaci\'on en Matem\'aticas IMUVA\\
E-mail: jesusjavier.jimenez@unican.es.
\\
\vskip.5cm
J.~Sanz:\\
Departamento de \'Algebra, An\'alisis Matem\'atico, Geometr{\'\i}a y Topolog{\'\i}a, Universidad de Valladolid\\
Facultad de Ciencias, Paseo de Bel\'en 7, 47011 Valladolid, Spain.\\
Instituto de Investigaci\'on en Matem\'aticas IMUVA\\
E-mails: jsanzg@am.uva.es.
\\
\vskip.5cm
G.~Schindl:\\
Fakult\"at f\"ur Mathematik, Universit\"at Wien,
Oskar-Morgenstern-Platz~1, A-1090 Wien, Austria.\\
E-mail: gerhard.schindl@univie.ac.at.

\end{document}